\documentclass[10pt]{amsart}

\usepackage[utf8]{inputenc}
\usepackage[english]{babel}
\usepackage{amsmath}
\usepackage{graphicx}
\usepackage{mathtools}
\usepackage{float}
\usepackage{amssymb}
\usepackage{esint}
\usepackage{amsfonts}
\usepackage[usenames, dvipsnames]{color}
\usepackage{multicol}
\usepackage{csquotes}
\usepackage{enumitem}
\usepackage{multicol}

\newtheorem{innercustomthm}{Theorem}
\newenvironment{customthm}[1]
  {\renewcommand\theinnercustomthm{#1}\innercustomthm}
  {\endinnercustomthm}
  
 \newtheorem{innercustomcorollary}{Corollary}
\newenvironment{customcorollary}[1]
  {\renewcommand\theinnercustomcorollary{#1}\innercustomcorollary}
  {\endinnercustomcorollary}

  \newtheorem{innercustomlemma}{Lemma}
\newenvironment{customlemma}[1]
  {\renewcommand\theinnercustomlemma{#1}\innercustomlemma}
  {\endinnercustomlemma}

\makeatletter
\DeclareFontFamily{U}{tipa}{}
\DeclareFontShape{U}{tipa}{m}{n}{<->tipa10}{}
\newcommand{\arc@char}{{\usefont{U}{tipa}{m}{n}\symbol{62}}}%

\newcommand{\arc}[1]{\mathpalette\arc@arc{#1}}

\newcommand{\arc@arc}[2]{%
  \sbox0{$\m@th#1#2$}%
  \vbox{
    \hbox{\resizebox{\wd0}{\height}{\arc@char}}
    \nointerlineskip
    \box0
  }%
}
\makeatother

  \newcommand{\deltah}{\delta^{\mathbb{H}}}
  
  \newcommand{\hp}{\phi^{\mathbb{H}}}
  \newcommand{\hsi}{\sigma^{\mathbb{H}}}
  
  \newcommand{\di}{\nabla \cdot}
  \newcommand{\hs}{\mathbb{H}^d_+}  
  \newcommand{\h}{\mathbb{H}}
  \newcommand{\bhs}{\partial \mathbb{H}^d_+}
  \newcommand{\rpb}{\partial \arc{B}^+}
  \newcommand{\fpb}{\partial \underline{B}^+}
  \newcommand{\Rp}{R^{\prime}}
  \newcommand{\exch}{\textrm{Exc}^{\mathbb{H}}}

\begin{document}

\title[Random Elliptic Operators on the Half-Space]{A Large-Scale Regularity Theory for Random Elliptic Operators on the Half-Space with Homogeneous Neumann Boundary Data}
\author{Claudia Raithel}
\begin{abstract}
In this note we derive large-scale regularity properties of solutions to second-order linear elliptic equations with random coefficients on the half-space with homogeneous Neumann boundary data; it is a companion to \cite{FischerRaithel} in which the situation for homogeneous Dirichlet boundary data was addressed.  Similarly to \cite{FischerRaithel}, the results in this contribution are expressed in terms of a first-order Liouville principle. It follows from an excess-decay that is shown through means of a stochastic homogenization-inspired Campanato iteration. The core of this contribution is the construction of a sublinear half-space-adapted corrector/vector potential pair that, in contrast to \cite{FischerRaithel}, is adapted to the Neumann boundary data.
\end{abstract}

\maketitle

\section{Introduction}
In this note we are interested in the large-scale boundary regularity of solutions to second-order linear elliptic equations with random coefficients and homogeneous Neumann boundary data. In particular, we work with the following model case: Let $u \in H^1_{loc}(\overline{\mathbb{H}}^d_+)$ be a weak solution of
\begin{align}
\label{intro1}
\begin{split}
-\di (a \nabla u) &= 0 \quad\quad\quad \textrm{ in }\quad \hs,
\\
e_d \cdot a \nabla u &= 0 \quad\quad\quad \textrm{ on }\quad \bhs,
\end{split}
\end{align}
where $a$ is the restriction to the half-space of a coefficient field $a(x): \mathbb{R}^d \rightarrow \mathbb{R}^{d \times d}$ that is bounded and uniformly elliptic on $\mathbb{R}^d$. This work is a continuation of \cite{FischerRaithel}, in which large-scale regularity properties for solutions to \eqref{intro1} with homogeneous Dirichlet boundary conditions were addressed. For simplicity we use scalar notation throughout, but our arguments also extend to the case of systems.

We recall that there are classical counterexamples showing that solutions $u \in H^1_{loc}(\mathbb{R}^d)$ of 
\begin{align}
 \label{introcorrect1}
-\di (a \nabla u) &= 0 \quad\quad\quad \textrm{ in }\quad \mathbb{R}^d,
\end{align}
where the coefficient field $a(x)$ is bounded and uniformly elliptic, may fail to be in $C_{loc}^{0, \alpha}(\mathbb{R}^d)$ for any $\alpha \in (0,1)$. In the scalar case this was shown, e.g., by Meyers \cite[Example 3]{PiccininiSpagnolo} and in the case of systems De Giorgi showed that solutions may even fail to be locally bounded \cite[Section 9.1.1]{GiaquintaMartinazzi}.

In contrast, if the coefficients $a$ in \eqref{introcorrect1} are spatially constant then solutions $u\in H^{1}_{loc}(\mathbb{R}^d)$ of \eqref{introcorrect1} are locally $C^{k, \alpha}$-H\"{o}lder continuous for all $k \in \mathbb{N}$ and $\alpha \in(0, 1)$. This can easily be deduced from the following two ingredients: a) the observation that $C^{1, \alpha}$-regularity is equivalent to a certain approximability of $u$ by linear functions (i.e. the equivalence of H\"{o}lder and Campanato spaces) and b) the decay of the tilt-excess given by
\begin{align}
\begin{split}
 \label{introcorrect2}
 & \fint_{B_{r}(x_0)} \left| \nabla u - \fint_{B_r(x_0)} \nabla u \, dx  \right|^2 \, dx\\
 & \quad \quad \quad \quad \quad \quad \quad \quad \quad \quad \leq C(d,\lambda)\left( \frac{r}{R}\right)^{2} \fint_{B_R(x_0)}\left| \nabla u - \fint_{B_R(x_0)} \nabla u \, dx \right|^2 \, dx, 
 \end{split}
\end{align}
which holds for all $ 0 < r \leq R$ and $x_0 \in \mathbb{R}^d$ (see, e.g., \cite[Theorem 5.14]{GiaquintaMartinazzi}).

Using the constant coefficient result, one may show that if $u$ solves \eqref{introcorrect1} with coefficients $a$ that are locally $C^{0,\alpha}$ then $u \in C_{loc}^{1,\alpha}(\mathbb{R}^d)$. To do this, one relies on a comparison of $u$ with the solution $v$ of \eqref{introcorrect1} on $B_R(x_0)$ with frozen coefficients; i.e., $v$ solves
\begin{align}
\label{introcorrect3}
\begin{split}
-\di (a(x_0) \nabla v) &= 0 \quad\quad\quad \textrm{ in }\quad B_R(x_0),
\\
v  &= u \quad\quad\quad \textrm{ on }\quad \partial B_R(x_0).
\end{split}
\end{align}
In particular, one applies \eqref{introcorrect2} to $v$ and uses the energy estimate for the equation satisfied by $u-v$ to show that there exists a ratio of radii $\theta= r/R$ such that the excess-decay
\begin{align}
 \begin{split}
  \label{campanato1}
  &\int_{B_r(x_0)}|\nabla u - \fint_{B_r}\nabla u \, dx |^2 \, dx \leq  C(d, \lambda, \| a \|_{C^{0,\alpha}(B_R(x_0))})\\
  & ~~~~~~~~~~~~~~~~~~~~ \times \left( \theta^{(2d+2+2\alpha)/2}   \int_{B_R(x_0)}|\nabla u - \fint_{B_R(x_0)}\nabla u \, dx|^2 \, dx +R^{d+2 \alpha} \right)
 \end{split}
\end{align}
holds for all $x_0 \in \mathbb{R}^d$. By iteration one obtains that 
\begin{align}
 \begin{split}
  \label{campanato1.1}
  &\int_{B_r(x_0)}|\nabla u - \fint_{B_r}\nabla u \, dx |^2 \, dx \leq C(d, \lambda, \| a \|_{C^{0,\alpha}(B_R(x_0))})\\
  & ~~~~~~~~~~~~~~~~~~~~  \times  \left(  \left( \frac{r}{R}\right)^{d+2\alpha}  \int_{B_R(x_0)}|\nabla u - \fint_{B_R(x_0)}\nabla u \, dx|^2 \, dx +r^{d+2 \alpha} \right)
 \end{split}
\end{align}
for all $0< r \leq R$. By observation a) in the previous paragraph one then obtains that $u \in C_{loc}^{1,\alpha}(\mathbb{R}^d)$. For general heterogeneous coefficient fields (with no continuity assumptions) this method-- a standard tool in regularity theory going by the name Campanato iteration (see, e.g., \cite[Lemma 5.13 and Theorem 5.19]{GiaquintaMartinazzi})-- fails as one may not view $u$ as a perturbation of $v$.

Recall that ``homogenization'' is said to occur if there exist constant coefficients $a_{hom}$ such that the operator $- \nabla \cdot a \nabla$ is well-approximated by $- \nabla \cdot a_{hom} \nabla$ in the macroscopic limit. When homogenization occurs then one may try to replicate the Campanato iteration described above in the case of H\"{o}lder continuous coefficients by instead viewing $u$ as a perturbation of the solution $v$ of the homogenized problem. Within the context of periodic homogenization this concept of a homogenization-inspired Campanato iteration has been around since the 80s and was originally used to obtain regularity results (up to the boundary) by Avellaneda and Lin in \cite{AvellanedaLinCPAM}. More recently there has been a large effort to extend this method to the case of random coefficient fields.

We now recall some basic notions from the homogenization theory of random linear elliptic operators. As one heuristically expects, homogenization occurs for coefficient fields with no long-range correlations that fluctuate at a scale much smaller than the macroscopic (material) scale. In stochastic homogenization one considers probability measures $\langle \cdot \rangle$ (called \textit{ensembles}) on the space of coefficient fields $a(x): \mathbb{R}^d \rightarrow \mathbb{R}^{d \times d}$ that are supported on the set
\begin{align}
\label{Omega}
\Omega & =\left. \left\{ a (x)  \right\vert   |a(x)\xi| \leq |\xi|\textrm{ and } \lambda |\xi|^2  \leq |\xi \cdot a(x) \xi| \textrm{ for all } \xi \in \mathbb{R}^d \textrm{ for a.e. } x \in \mathbb{R}^d\right\}
\end{align}
and for which homogenization occurs almost surely. In particular, one assumes that the measure $\langle \cdot \rangle$ satisfies the following two properties:
\begin{itemize}
 \item \textit{Stationarity}: The measure $\langle \cdot \rangle$ must be shift invariant, where for $x \in \mathbb{R}^d$ the shift $\tau_x: \Omega \rightarrow \Omega$ is given by $\tau_x(a(\cdot)) \mapsto a(\cdot + x)$. 
 \item \textit{Ergodicity}: Shift invariant random variables must $\langle \cdot \rangle$-almost surely be constant. Morally, this corresponds to the qualitative decorrelation of the coefficient field on large-scales. As the coefficient fields constructed by Meyers and De Giorgi in their counterexamples are radial they are in this sense ``ungeneric''. In this note, we require a slightly quantified version of ergodicity (see condition \eqref{reqwsc}) for our arguments to work. This requirement is the same as the one used in \cite{FischerRaithel}.
\end{itemize}

A central object in homogenization is the \textit{corrector}. It ``corrects" the linear coordinate functions $x \in \mathbb{R}^d  \mapsto x_i$ as to be $a$-harmonic; so, the whole-space corrector in the direction $\xi \in \mathbb{R}^d$, denoted $\phi_{\xi}$, is a distributional solution of
\begin{align}
 \label{wsceqp}
 - \di( a \nabla (\phi_{\xi} + \xi \cdot x)) = 0 \quad \quad \quad \textrm{ in} \quad \mathbb{R}^d.
\end{align}
In a certain sense, the corrector ``cancels out" the fluctuations of the heterogeneous coefficients and, therefore, oscillates at the same scale. 

Clearly, solutions of (\ref{wsceqp}) are not unique. To determine a unique (up to addition of a constant) choice of corrector, we ask that in the expression $\phi_{\xi} + \xi \cdot x$ the linear function should be dominant on large scales. In particular, we are interested in solutions of \eqref{wsceqp} that are \textit{sublinear}. The almost sure existence and uniqueness (up to the addition of a constant) of a sublinear corrector is established by Gloria, Neukamm, and Otto in \cite{GloriaNeukammOtto} assuming only stationarity and qualitative ergodicity. This corrector is actually shown to be jointly sublinear with a corresponding vector potential $\sigma$ (see \eqref{wsceqs}) in the sense that the averages over balls $\delta_{r}^{GNO}(\phi, \sigma)$, given by
\begin{align}
\label{defndeltawscGNO}
\begin{split}
&\delta^{GNO}_{r}(\phi, \sigma)=  \\
&~~~\frac{1}{r} \left( \fint_{B_r}\sum_{i=1}^d \left(|\phi_{e_i} -\fint_{B_r} \phi_{e_i} \, dx  \,|^2   + \sum_{j,k = 1}^d |\sigma_{e_ijk} -\fint_{B_r} \sigma_{e_ijk } \, dx |^2  \right) \, dx \right)^{1/2},
\end{split}
\end{align}
satisfy the condition $\lim_{r\rightarrow \infty}\delta^{GNO}_{r}(\phi, \sigma)=0$. While their construction of the corrector is similar to previous treatments (see, e.g., \cite[Section 7.2]{Homogbook}), \cite{GloriaNeukammOtto} is the first instance of the vector potential $\sigma$, which is a common object in periodic homogenization, being used in setting of random coefficients.

The construction of the corrector in \cite{GloriaNeukammOtto} also guarantees that the random field $\nabla \phi$ is stationary. This allows one to use a heuristic argument to obtain a relation for the homogenized coefficients. In particular, we recall that for an ergodic ensemble the spatial average of a stationary random field may be replaced by taking the expectation at a fixed point. As the homogenized macroscopic current $a_{hom} \xi$ must coincide with the spatial average of the corrected microscopic current $a(\xi + \nabla \phi_{\xi})$, this observation yields that
\begin{align}
\label{homogenizedcoeff}
a_{hom} \xi =  \mathbb{E}\left[a(\xi + \nabla \phi_{\xi})\right].
\end{align} 
The \textit{current correction} in the direction $\xi$ is then given by $q_{\xi} = a (\nabla \phi_{\xi}+ \xi) - a_{hom} \xi$.

We may now define the vector potential of the current correction $\sigma_{\xi}$, which is taken to be a skew-symmetric (in $j$ and $k$) distributional solution of 
\begin{align}
 \label{wsceqs}
 \nabla_k \cdot \sigma_{\xi jk} = q_{\xi j} \quad \quad \quad \textrm{in} \quad \mathbb{R}^d.
\end{align}
The vector potential as defined above is only unique up to the addition of solenoidal fields. However, in \cite{GloriaNeukammOtto} it is shown that by choosing the correct gauge one obtains a unique $\sigma_{\xi}$  (up to addition of a constant) such that the joint sublinearity condition with $\phi_{\xi}$ is satisfied. 

It is standard in homogenization to, on any ball $B_R$, approximate the $a$-harmonic function $u$ with the solution $u_{hom}$ of the homogenized problem, which solves
\begin{align}
\label{homogenized}
\begin{split}
-\di (a_{hom} \nabla u_{hom}) &= 0 \quad\quad\quad \textrm{ in }\quad B_R,
\\
u_{hom}& = u \quad\quad\quad \textrm{ on }\quad \partial B_R,
\end{split}
\end{align}
using a two-scale asymptotic expansion that is truncated to yield
\begin{align}
\label{twoscale}
u \approx u_{hom}+  \sum_{i=1}^d \phi_{e_i} \partial_{x_i} u_{hom}.
\end{align}
The two scales that appear here are the macroscopic scale at which $u_{hom}$ changes (the ``slow'' scale) and the microscopic scale at which the heterogeneous coefficients oscillate (the ``fast'' scale). The error in this approximation, which we denote as $w = u - (u_{hom} + \sum_{i=1}^d\phi_{e_i}\partial_{x_i} u_{hom})$, solves 
\begin{align}
 \label{divformerror}
 - \nabla \cdot (a \nabla w) = \nabla \cdot \left( \sum_{i=1}^d (\phi_{e_i}a - \sigma_{e_i})\partial_{x_i} \nabla u_{hom}\right)\quad\quad\quad\textrm{in}\quad B_R.
\end{align}

In terms of large-scale regularity results, in this note we emphasize Liouville principles. Notice that if an excess-decay of the type \eqref{introcorrect2} holds on large scales (i.e. when $r^* \leq r \leq R$ for some $r^*>0$) then this is sufficient to prove a first-order Liouville principle: One may show that the dimension of the space of $a$-harmonic functions satisfying the growth condition $|u(x)| \leq C(1 + |x|^{1+\alpha})$ for some $C \in \mathbb{R}$ and $\alpha \in (0,1)$  is the same as in the Euclidean setting (when $a=\textrm{Id}$). The link between regularity results and Liouville statements is further seen in that the counterexamples of Meyers and De Giorgi are both also counterexamples to a zeroth-order Liouville principle (they are sublinear and $a$-harmonic, but not constant).

As previously mentioned, improving large-scale regularity results for random linear elliptic operators using homogenization results has been an active area and we now review some literature. The concept was first seen in \cite{BenjaminiCopinKozmaYadin} within the context of random walks in random environments; In this work, Benjamini, Duminil-Copin, Kozma, and Yadin prove a zeroth-order Liouville principle in the setting of a supercritical percolation cluster under the assumption of qualitative ergodicity. Shortly afterwards, Marahrens and Otto, in the more analytic contribution \cite{MarahrensOtto}, obtained a large-scale $C^{0, \alpha}$-regularity theory for $\alpha \in (0,1)$ assuming an ensemble that satisfies a logarithmic Sobolev inequality. A large-scale  $C^{0,1}$-theory was then obtained for scalar equations under a finite range of dependence assumption by Armstrong and Smart in \cite{ArmstrongSmart}. The work by Armstrong and Smart was the first in which the scheme of Avellaneda and Lin from \cite{AvellanedaLinCPAM} was adapted to the setting of random coefficients. The contribution \cite{ArmstrongSmart} was then followed by \cite{ArmstrongMourrat}, in which Armstrong and Mourrat were able to replace the finite range of dependence condition by a weaker ``$\alpha$-mixing'' condition and also treat the case of systems. The contribution \cite{ArmstrongMourrat} was predated by the first version of \cite{GloriaNeukammOtto}, which motivated both the present contribution and \cite{FischerRaithel}. In \cite{GloriaNeukammOtto} Gloria, Neukamm, and Otto obtained a large-scale $C^{1,\alpha}$-regularity theory for $\alpha \in (0,1)$ and a first-order Liouville principle under only a qualitative ergodicity assumption. The regularity results in \cite{GloriaNeukammOtto} are obtained through means of a Campanato iteration that hinges on the existence of a sublinear corrector/vector potential pair.

Using a slight quantification of ergodicity -- namely a growth condition on the corrector-- Fischer and Otto were then able to obtain a large-scale $C^{k, \alpha}$-theory in \cite{FischerOtto}. Their quantitative assumption on the sublinear growth of the corrector is required for the construction of higher-order correctors. The construction of the higher-order correctors in \cite{FischerOtto} inspired the construction of the half-space-adapted corrector in \cite{FischerRaithel} and also in this note. It should be noted that \cite{FischerOtto} was followed by \cite{ArmstrongKuusiMourratComplete} in which a different proof of a large-scale $C^{k, \alpha}$-regularity theory is given; the results of Armstrong, Kuusi, and Mourrat in \cite{ArmstrongKuusiMourratComplete} are valid for ensembles satisfying the $\alpha$-mixing condition of \cite{ArmstrongMourrat}. Also, recently, Armstrong and Dario extended the results of \cite{BenjaminiCopinKozmaYadin} and obtained a large-scale $C^{k, \alpha}$-regularity theory on supercritical percolation clusters; In particular, they proved higher-order Liouville principles.

While to the best of our knowledge \cite{FischerRaithel} was the first instance of a half-space corrector being constructed within the setting of random coefficients, in \textit{periodic} homogenization boundary correctors for the Dirichlet problem were already introduced by Avellaneda and Lin in \cite{AvellanedaLinCPAM}. Twenty years later, Kenig, Lin, and Shen introduced  boundary correctors for the Neumann problem in \cite{KenigLinShenNeumann} and were able to extend the uniform (in $\varepsilon$) Lipschitz and $W^{1,p}$ estimates for solutions $u_{\varepsilon}$ of the Dirichlet problem
\begin{align}
\label{uepsilon}
\begin{split}
-\di (a\left(\frac{x}{\varepsilon}\right) \nabla u^{\varepsilon}) &= f \quad\quad\quad \textrm{ in }\quad \Omega,
\\
u^{\varepsilon} &= g \quad\quad\quad \textrm{ on }\quad \partial \Omega,
\end{split}
\end{align}
shown by Avellaneda and Lin in \cite{AvellanedaLinCPAM} to the Neumann setting. One thing which makes the Neumann situation more complicated than the Dirichlet case is that the no-flux boundary condition for $u^{\varepsilon}$ is actually $\varepsilon$-dependent. In the subsequent work \cite{KenigLinShenGreen}, Kenig, Lin, and Shen then used the uniform estimates and techniques developed in \cite{KenigLinShenNeumann} and the results of \cite{AvellanedaLinCPAM} to study the asymptotics (as $\varepsilon \rightarrow$ 0) of the Green and Neumann functions. This allowed them to obtain near optimal first-order convergence rates for $u^{\varepsilon}\rightarrow u_{hom}$ as $\varepsilon \rightarrow 0$ for both the Dirichlet and Neumann problems in $W^{1,q}$ for $1 <q\leq\infty$ and also to give a more refined estimate for the error in the asymptotic expansion of the Poisson kernel given in \cite{AvellanedaLinJMPA}. It should also be mentioned that in \cite{ArmstrongShen} Armstrong and Shen obtained a $C^{0,1}$-regularity theory up to the boundary for both the Dirichlet and Neumann problems in the \textit{almost periodic} setting using arguments that likely extend to the case of random coefficients. 

We would also like to mention that in the homogenization of linear elliptic equations on bounded domains, as already alluded to above, one runs into a boundary layer phenomenon.  Returning to the setting of \textit{periodic} coefficients, in \cite{KenigLinShenGreen} the authors replace the standard correctors by the boundary correctors in the 2-scale expansion to handle this boundary layer and, thereby, obtain their improved convergence rates (the classical estimate is $O(\varepsilon^{1/2})$, whereas one would expect and Keing, Lin, and Shen almost obtain $O(\varepsilon)$). As is already noted in \cite{KenigLinShenGreen}, their results on the asymptotic behaviour of the Poisson kernel can be used to investigate the \textit{oscillating Dirichlet problem}, in which not only the coefficients oscillate at a scale $\varepsilon$, but also the boundary data.  The desire to consider the homogenization of the oscillating Dirichlet problem stems from attempting to derive higher-order convergence rates for the homogenization of the standard Dirichlet problem.  Obtaining results concerning the homogenization of the oscillating Dirichlet problem is, however, much more difficult than in the standard Dirichlet case. In particular, the situation turns out to depend on the geometry of the domain; Specifically, on whether the tangent hyperplanes are resonant with the periodic structure of the coefficients.

There is a large body of literature concerned with higher-order convergence rates for the homogenization of the standard Dirichlet problem in \textit{periodic} homogenization. We would only like to mention a couple of works, the first of which is \cite{AllaireAmar}. Here, Allaire and Amar treat the special case of a $\mathbb{Z}^d$-periodic coefficient field, where the domain is taken to be the open unit cube $(0,1)^d$ and $\varepsilon^{-1}\in \mathbb{N}$. Their work relies on the introduction of a first-order boundary layer term $u_{1}^{bl, \varepsilon}(x)$, which basically corrects the first-order 2-scale expansion such that the ansatz given for $u_{\varepsilon}$ has the appropriate boundary data. The boundary layer term is characterized as the solution of an \textit{oscillating Dirichlet problem}, which Allaire and Amar then treat by splitting it into two contributions: one which may be treated with the homogenization results available for the standard Dirichlet problem and another that decays to $0$ in the interior of the domain as $\varepsilon \rightarrow 0$.  Ultimately, they are able to give functions $u_1(x, x /\varepsilon)$ and $u_2(x, x/ \varepsilon)$ such that for any $\omega \subset \subset \Omega$ the rate $\|u_{\varepsilon} - u_{hom} - \varepsilon u_1(x, x/\varepsilon) - \varepsilon^2 u_2(x, x/\varepsilon) \|_{H^1(\omega)} = O( \varepsilon^{3/2})$ holds.


The work of Allaire and Amar was followed by \cite{GerardVaretMasmoudi2} in which G\'{e}rard-Varet and Masmoudi examined the case of polygonal domains with normals satisfying a diophantine condition. Their strategy is to approximate $u_1^{bl, \varepsilon}$ in terms of functions that are $a$-harmonic on the half-planes intersected to obtain the domain with prescribed $\mathbb{Z}^d$-periodic boundary data. In the case of a polygon with sides of rational slope, for each such half-space problem the periodicity is retained in the directions tangential to the half-plane; Of course, when the slopes are diophantine this is no longer true. In \cite{GerardVaretMasmoudi2} the authors use their diophantine assumption to treat the half-space problems and then glue the solutions together to approximate $u_{1}^{bl,\varepsilon}$. In this way, they are able to find functions $u_1(x,x/ \varepsilon)$ and $u_2(x, x / \varepsilon)$ such that for all $\omega \subset \subset \Omega$ the rate $\|u_{\varepsilon} - u^0  - \varepsilon u^1(x, x / \varepsilon) - \varepsilon^2 u^2(x, x / \varepsilon)\|_{H^1(\omega)}  = O(\varepsilon^2)$ holds.  In \cite{GerardVaretMasmoudi} the same authors then considered uniformly convex domains and obtained the rate $\|u_{\varepsilon} - u_{hom} - \varepsilon u_1(x, x/\varepsilon) - \varepsilon^2 u_2(x, x/\varepsilon) \|_{H^1(\omega)} = O( \varepsilon^{1 + \alpha})$ for $\alpha < \frac{d-1}{3d+5}$. In both \cite{GerardVaretMasmoudi2} and \cite{GerardVaretMasmoudi}, quantitative proofs are given for the homogenization of the oscillating Dirichlet problem.

To conclude our discussion of the boundary layer in \textit{periodic} homogenization, we mention a couple further works concerned with convergence rates for the oscillating Dirichlet problem in uniformly convex domains. We first remark that in the case of constant coefficients and oscillating boundary data,  Aleksanyan, Shahgholian, and Sjolin were able to obtain $L^q$- convergence rates for $2 \leq q <\infty$ and $d \geq 2$ in \cite{AleksanyanShahgholianSjolin}; For dimensions $d \geq 4$ they obtain the optimistically hypothesised rate $O(\varepsilon)$. More recently, the story was basically completed in \cite{ArmstrongKuusiMourratPrange} in which Armstrong, Kuusi, Mourrat, and Prange obtained $L^q$-convergence rates for $2 \leq q <\infty$ and $d \geq 2$ that are nearly optimal for $d \geq4$ in that the agree with the rates of \cite{AleksanyanShahgholianSjolin} up to the loss of an arbitrarily small exponent $\delta>0$. Their analysis makes use of the two-scale expansion for the Poisson kernel proved by Kenig, Shen, and Lin in \cite{KenigLinShenGreen} and relies on the treatment of the half-space problems mentioned in the context of \cite{GerardVaretMasmoudi} above. In their argument they, in particular, approximate the Dirichlet corrector used by Kenig, Lin, and Shen in terms of cell-correctors and half-space boundary layer correctors using the Lipschitz theory of Avellaneda and Lin in \cite{AvellanedaLinCPAM}. As is mentioned in \cite{ArmstrongKuusiMourratPrange}, it was noticed by Shen and Zhuge in \cite{ShenZhuge} that, after upgrading the regularity result obtained in \cite{ArmstrongKuusiMourratPrange} for the homogenized boundary data, the strategy of \cite{ArmstrongKuusiMourratPrange} also yields nearly optimal rates for the dimensions $d=2,3$.  Lastly, we remark that, in contrast to the Dirichlet case, the oscillating Neumann problem with zero-order oscillating boundary data (i.e. $n \cdot a(x/\varepsilon) \nabla u_{\varepsilon}  = g_0(x, x/ \varepsilon)$ on $\partial \Omega$) has been well-understood for sometime (see, e.g., \cite[Chapter 1, Section 7.1]{BensoussanLionsPapanicolaou}). Also in \cite{ShenZhuge}, towards obtaining higher-order convergence rates for the homogenization of the standard Neuman problem, the case of first-order oscillating boundary data  (i.e. $n \cdot a(x/\varepsilon) \nabla u_{\varepsilon}  = T_{ij} \cdot \nabla_x (g_{ij}(x, x/ \varepsilon)) + g_0(x, x/ \varepsilon)$ on $\partial \Omega$ with $T_{ij} = n_ie_j - n_j e_i$) is treated. Following an approach inspired by \cite{ArmstrongKuusiMourratPrange}, Shen and Zhuge obtain the nearly optimal $L^q$-convergence rates for $2 \leq q <\infty$ and $d \geq 3$.

Returning to the current work: In this contribution we derive large-scale regularity results analogous to those in \cite{GloriaNeukammOtto} for the whole-space case and in \cite{FischerRaithel} for the half-space case with homogeneous Dirichlet boundary data. We prove two theorems: In the first theorem, assuming that for a given realization $a$ there exists a whole-space corrector/vector potential pair with a certain quantitative sublinear growth property, we construct a sublinear half-space-adapted corrector/vector potential pair, which we denote $(\hp, \hsi)$. Then, in the second theorem,  we use the sublinearity of $(\hp, \hsi)$ to prove a large-scale intrinsic $C^{1, \alpha}$ excess-decay for solutions of \eqref{intro1} with $\alpha \in (0,1)$ and the tilt-excess adapted to both the homogenization and half-space settings.

To determine the appropriate notion of tilt-excess for our setting we use a heuristic observation that arises due to the presence of the no-flux boundary condition in \eqref{intro1}. Let
\begin{align}
\label{defnB}
B := \left\{b \in \mathbb{R}^d \hspace{.2cm} | \hspace{.2cm} e_d \cdot a_{hom} b = 0 \right\}
\end{align}
\noindent and denote an orthonormal basis of $B$ as  $\left\{b_1,..., b_{d-1} \right\}$. We may complete this to an orthonormal basis of $\mathbb{R}^d$ with some $b_d \in \mathbb{R}^d$. Due to the boundary condition in \eqref{intro1} we find that it suffices to compare $u$ to the space of $a$-affine functions without a component in $b_d$-direction. This motivates the definition
\begin{align}
\label{defnexcess}
\textrm{Exc}^{\h}(r):=\inf_{ \tilde{b} \in B} \fint_{B_r^+} |\nabla u-( \tilde{b} + \nabla \hp_{\tilde{b}})  |^2 \,dx.
\end{align}
For a further motivation of \eqref{defnexcess} recall that we seek to prove a first-order Liouville principle for solutions of \eqref{intro1}. The Liouville principle will, of course, have to hold in the constant coefficient case when $a = a_{hom}$ and subquadratic solutions of \eqref{intro1} are of the form $\tilde{b} \cdot x + c$ for $\tilde{b} \in B$ and $c \in \mathbb{R}$. Since the excess should compare $u$ to the space that we expect to characterize the subquadratic solutions of \eqref{intro1} for general $a$, we again arrive at \eqref{defnexcess}. This definition hints that it will only be necessary to construct the half-space-adapted corrector/ vector potential pair in the directions $b_i$ for $i \neq d$.

As already mentioned above, in order to obtain our results for a given realization of the random coefficient field it suffices to assume the existence of a whole-space corrector/vector potential pair satisfying a slightly quantified sublinearity condition. In particular, we assume that
\begin{align}
\label{reqwscGNO}
 \sum_{m=0}^{\infty} m \cdot  \delta^{GNO}_{2^m}(\phi, \sigma)^{1/3} < \infty.
\end{align}
It is simple to see that this assumption actually implies that 
\begin{align}
\label{reqwsc}
 \sum_{m=0}^{\infty} m \cdot  \delta_{2^m}(\phi, \sigma)^{1/3} < \infty,
\end{align}
where
\begin{align}
\label{defndeltawsc}
\delta_{r}(\phi, \sigma)= \frac{1}{r} \left( \fint_{B_r}\sum_{i=1}^d \left(|\phi_{e_i} \,|^2   + \sum_{j,k = 1}^d |\sigma_{e_ijk} |^2  \right) \, dx \right)^{1/2}.
\end{align}
As both of the maps $\xi \mapsto \phi_{\xi}$ and $\xi \mapsto \sigma_{\xi}$ are linear, using the definition of $\delta_{r}(\phi, \sigma)$ given in \eqref{defndeltawsc} and Young's inequality, we find that for any orthonormal basis $\left\{b_1,..., b_d \right\}$ of $\mathbb{R}^d$  
\begin{align}
\label{different_basis}
\begin{split}
&\frac{1}{r} \left( \fint_{B_r}\sum_{i=1}^d \left(|\phi_{b_i} \,|^2   + \sum_{j,k = 1}^d |\sigma_{b_ijk} |^2  \right) \, dx \right)^{1/2}\\
 = & \frac{1}{r} \left(  \fint_{B_r}\sum_{i=1}^d \left(|\sum_{w=1}^d \langle b_i, e_w\rangle \phi_{e_w} \,|^2   + \sum_{j,k = 1}^d |\sum_{w=1}^d \langle b_i, e_w \rangle \sigma_{e_wjk}|^2  \right) \, dx \right)^{1/2}\\
 \leq & \frac{1}{r} \left( \frac{d+1}{2} \fint_{B_r}\sum_{i=1}^d \left(\sum_{w=1}^d |\langle b_i, e_w\rangle|^2 |\phi_{e_w} \,|^2   + \sum_{j,k,w = 1}^d   |\langle b_i, e_w \rangle|^2 | \sigma_{e_wjk}  |^2  \right) \, dx \right)^{1/2}\\
 \leq &\left( \frac{d+1}{2}  \right)^{1/2} \frac{1}{r} \left(  \fint_{B_r}\sum_{i=1}^d  \left( \sum_{w=1}^d |\phi_{e_w} \,|^2   + \sum_{j,k,w = 1}^d | \sigma_{e_wjk}  |^2  \right) \, dx \right)^{1/2}\\
 \leq & \left( \frac{d(d+1)}{2}  \right)^{1/2} \delta_{r}(\phi, \sigma)\\
 \end{split}
\end{align}
holds with $\langle \cdot, \cdot \rangle $ denoting the standard dot product. 

As was already noted in \cite{FischerRaithel}, it can be shown that (\ref{reqwsc}) is satisfied $\langle \cdot \rangle$- almost surely for a large class of stationary and ergodic ensembles; for example, when the coefficient field $a(x)$ has a finite range of dependence (that (\ref{reqwsc}) is satisfied $\langle \cdot \rangle$- almost surely follows from \cite{GloriaOttoNew}) or is the image under a Lipschitz mapping $\psi: \mathbb{R}^{d \times d} \rightarrow \Omega$ of a matrix-valued stationary Gaussian random field whose spatial correlations satisfy a prescribed slow decay (that (\ref{reqwsc}) is satisfied $\langle \cdot \rangle$- almost surely follows from  \cite{CorrectorEstimatesSlowDecorrelation}).\\

\noindent \textbf{Remark}: The general layout and strategy in this note resemble \cite{FischerRaithel}. However, there are some differences; the most prominent of these may be found in Steps 1 and 2 of Section 3.\\

\noindent\textbf{Notation}: By $ ``a \lesssim b'' $ we always mean $`` a \leq C(d, \lambda) b ''$, where $C(d, \lambda)$ is an arbitrary constant depending on the dimension $d$ and the ellipticity ratio $\lambda$. We call the $i$-th coordinate vector $e_i$ so that $(e_i)_j = \delta_{ij}$. Furthermore, we use the Einstein summation convention under which an index is summed over if it appears twice. For example, using this convention, by $\langle b_i, \nabla v \rangle b_i$ we mean $\sum_{i=1}^d \langle b_i, \nabla v \rangle b_i $. Occasionally, we may include the summation symbol to avoid confusion.

For a measurable set $ V \subseteq \mathbb{R}^d$ we use $\chi_{V}$ to denote the indicator function. The Lebesgue measure of $V$ is written as $|V|$. We use $C_c^{\infty}(V)$ to denote smooth functions with compact support in $V$. By $u \in H^1_{loc}(\overline{\mathbb{H}}^d_+)$ we mean that $u \in H^1(V\cap \mathbb{H}^d_+)$ for any bounded, open set $V \subseteq \mathbb{R}^d$. We, furthermore, use the notation $H^1_{bdd}(\hs) = \left\{ u \in H^1(\hs): \textrm{supp}(u)\subseteq B_r \textrm{ for some } r>0 \right\}$. 

We denote $B_r:= \{x\in \mathbb{R}^d \, | \, |x| < r \}$ and $B_r^+ := \left\{ x\in \mathbb{R}^d \, | \, |x| < r  \textrm{ and } x_d>0 \right\}$. The boundary of the half-ball $B_r^+$ is decomposed into round and flat parts: $\rpb_r = \partial B_r \cap \hs $ is the round part and $\fpb_r = B_r \cap \partial \hs $ is the flat part.\\

\section{Main Results}
We now give the full statement of the two theorems and the Liouville principle that arises as a corollary. In this first theorem we construct the half-space-adapted corrector/vector potential pair:

\begin{customthm}{1}\label{existenceofcorrectors}
Let $a \in \Omega$ and let $\left\{b_i \right\}$ be an orthonormal basis of $\mathbb{R}^d$ such that $b_i \in B$ for $i \neq d$. Here, $B$ is given by \eqref{defnB} and $\Omega$ is defined in \eqref{Omega}. Assume that there exists a whole-space corrector/vector potential pair $(\phi, \sigma)$ satisfying the equations (\ref{wsceqp}) and (\ref{wsceqs}) along with the additional growth condition \eqref{reqwsc}. Then there exists a half-space-adapted corrector/vector potential pair $(\hp, \hsi)$ satisfying the following properties:

\begin{itemize}
    \item[i)] The half-space-adapted corrector $\hp_{b_d}$ and the half-space-adapted vector potential $\hsi_{b_d}$ are the restriction of $\phi_{b_d}$ and $\sigma_{b_d}$ respectively to the half-space, i.e. $\hp_{b_d} = \phi_{b_d}|_{\hs}$ 
    and $ \hsi_{b_d} = \sigma_{b_d}|_{\hs}$.
    \item[ii)] For $i \neq d$ the half-space-adapted corrector $\hp_{b_i}$ is a weak solution of
   \begin{subequations}
    \label{halfspacecorrector}
    \begin{align}
    \label{CorrectorEquationHalfSpace}
    ~~~~~~~~-\di a \nabla (\hp_{b_i} + b_i \cdot x) &= 0 &&\textrm{in }\quad \hs,
    \\
    \label{boundaryconditions}
    ~~~~~~~~ e_d \cdot a \nabla (\hp_{b_i} + b_i \cdot x) & =  0 &&\textrm{on }\quad \bhs
    \end{align}
    \end{subequations}
    where the class of test functions is given by $H^1_{bdd}(\hs)$.
    \item[iii)] For $i \neq d$ and $j \in \left\{1,...,d \right\}$ the half-space-adapted vector potential $\hsi_{b_ij}$ is a distributional solution of
    \begin{align}
    \label{halfspacepotential}
    \nabla_k \cdot \hsi_{b_ijk} = e_j \cdot \left(a(\nabla \hp_{b_i} + b_i) - a_{hom} b_i \right) \quad \quad \quad \textrm{in} \quad \hs.
    \end{align}
     Furthermore, $\hsi_{b_i jk}$ is skew-symmetric in $j$ and $k$.
    \item[iv)] The half-space-adapted corrector/vector potential pair $(\hp, \hsi)$ is sublinear in the sense that
    \begin{align}
    \label{defndeltah}
    \begin{split}
    & \deltah_{r}(\hp, \hsi) := \frac{1}{r} \left( \sum_{i=1}^{d-1}  \fint_{B_r^+} |(\hp_{b_i} - \fint_{B_r^+}\hp_{b_i}\, dx, \hsi_{b_i})|^2 \,dx \right. \\
    & \left.\quad \quad \quad \quad \quad \quad \quad \quad \quad \quad \quad \quad \quad \quad \quad \quad + \fint_{B_r} |(\phi^{\h}_{b_d}, \sigma^{\h}_{b,d})|^2 \,dx \vphantom{\sum_{i=1}^{d-1}  \fint_{B_r^+} |(\hp_{b_i} - \fint_{B_r^+}\hp_{b_i}\, dx, \hsi_{b_i})|^2 \,dx}\right)^{1/2}
    \end{split}
    \end{align}
    satisfies
    \begin{align}
    \label{sublinear}
    \lim_{r\rightarrow \infty} \deltah_r (\hp, \hsi) =0.
    \end{align}
\end{itemize}
\end{customthm}

\vspace{.2cm}

Actually, the estimates used to prove Theorem 1 guarantee a growth rate for the half-space-adapted corrector/vector potential pair. Just like in \cite{FischerRaithel}, if the whole-space corrector/ vector potential pair is sublinear in the sense that it satisfies
\begin{align*}
 \delta_r \lesssim \frac{1}{r^{\gamma}}
\end{align*}
for $\gamma >0$ then the estimates \eqref{rate1} and \eqref{pfthm13} ensure that the half-space-adapted corrector/vector potential pair satisfies
\begin{align*}
 \deltah_r \lesssim \frac{1}{r^{\gamma/3}}.
\end{align*}

The sublinear pair $(\hp, \hsi)$ constructed in Theorem 1 is then used to prove Theorem 2.

\begin{customthm}{2}\label{excessdecay}  Let $a \in \Omega$. Then for all H\"{o}lder exponents $\alpha \in (0,1)$ there exists a constant $C_{\alpha}(d, \lambda)$ such that if for a radius $R>0$ there exists a half-space-adapted
corrector/vector potential pair satisfying $i)-iii)$ from Theorem 1 on $B_R^+$ and there exists a minimal radius $r^*_{\alpha}<R$ for which
\begin{align}
 \label{smallness} 
 \deltah_r(\hp, \hsi) \leq \frac{1}{C_{\alpha}(d, \lambda)} \textrm{ if } r>r^*_{\alpha},
\end{align}
\noindent the following properties hold:\\

Let $u\in H^1(B_R^+)$ be $a$-harmonic with no-flux boundary conditions on $\fpb_{R}$, i.\,e.\ let $u$ be a weak solution of
\begin{align*}
-\di (a \nabla u) &= 0 \quad\quad\quad \textrm{ in }\quad B_R^+,
\\
e_d \cdot a \nabla u &= 0 \quad\quad\quad \textrm{ on }\quad \fpb_R
\end{align*}
\noindent where the class of test functions is given by $\left\{ u \in H^1(\hs) \, : \, \textrm{supp}(u) \subset B_r \textrm{ for some } \right.$ $\left.  R> r>0 \right\}$. We define the half-space-adapted tilt-excess of $u$ on the half-ball $B_r^+$ as indicated in \eqref{defnexcess}.\\

\noindent Then:

\begin{itemize}
 \item [i)]For $r\in [r^\ast_{\alpha},R]$  the excess-decay estimate given by
\begin{align}
\label{ExcessDecay}
 \textrm{Exc}^{\h}(r) \lesssim \left( \displaystyle\frac{r}{R} \right)^{2\alpha} \textrm{Exc}^{\h}(R)
\end{align}
\noindent holds.

 \item [ii)]For $r \in [r^{\ast}_{\alpha}, R]$ the tilt-excess functional
\begin{align*}
\tilde{b} \in \mathbb{R}^d \mapsto \fint_{B_r^+} |\nabla u -( \tilde{b} + \nabla \phi_{\tilde{b}}) |^2 \,dx
\end{align*}
is  coercive.

\item[iii)] There exists $C_{Mean}(d, \lambda)\geq1$ such that for $ r \in [ r^\ast_{1/2}, R]$ the mean-value property 
\begin{align}
 \label{mvp}
 \fint_{B_r^+} |\nabla u|^2 \,dx \leq  C_{Mean} \fint_{B_R^+}| \nabla u|^2 \,dx
\end{align} 
holds.
\end{itemize}
\end{customthm}

We obtain the following Liouville statement as a corollary:

\begin{customcorollary}{1}\label{Liouville}
For $a \in \Omega$ assume that there exists a whole-space corrector/vector potential pair $(\phi, \sigma)$ satisfying the growth condition (\ref{reqwsc}) and let $B$ be given by \eqref{defnB}. Then the following first-degree Liouville principle holds:\\

If  $u\in H^1_{loc}(\overline{\mathbb{H}}^d_+)$ is an $a$-harmonic function with no-flux boundary conditions on $\bhs$ (i.e. $u$ is a solution of \eqref{intro1}), which grows subquadratically in the sense that 
\begin{align}
\label{subquad}
\lim_{r \rightarrow \infty}\frac{1}{r^{1+\alpha}}\bigg(\fint_{B_r^+}{|u|^2} \,dx \bigg)^{1/2}=0 \textrm{ for some } \alpha>0
\end{align}
\noindent then $u$ is of the form $u = \tilde{b} \cdot x + \hp_{\tilde{b}} + c$ for some $\tilde{b} \in B$ and $c \in \mathbb{R}$.
\end{customcorollary}

In particular, for ensembles guaranteeing the $\langle \cdot \rangle$-almost sure existence of a whole-space pair $(\phi, \sigma)$ satisfying the growth condition \eqref{reqwsc} this gives a $\langle \cdot \rangle$- almost sure first-order Liouville principle.

\section{construction of the half-space-adapted corrector/vector potential pair}

\noindent \textbf{Notation}: We denote the homogeneous Sobolev space of functions with a square-integrable gradient on the half-space as $\dot{H}^1(\hs)$. We, furthermore, define $\dot{H}_0^1(\hs) = \{ u \in \dot{H}^1(\hs) \hspace{.2cm} | \hspace{.2cm} u = 0 \textrm{ on } \bhs \}$ equipped with the inner-product $\langle u, v \rangle_{\dot{H}^1} =  \int_{\hs} \nabla u \cdot \nabla v \, dx$. 

We never use a subscript of the form ``$,N$'' to refer to a partial derivative. Instead, the subscripts that appear after a comma refer to a scale; e.g., $\hp_{b_i, N}$ is an intermediate half-space-adapted corrector that has the desired boundary condition \eqref{boundaryconditions} on $\fpb_{r_02^{N}}$ (where $r_0$ is an initial radius to be chosen later).

We must sometimes keep track of explicit constants. Therefore, we introduce the following definitions:
\begin{itemize}[topsep=0pt]
\item $C_P(d)$ denotes the Poincar\'{e} constant on $B_1^+$ for functions with zero average.
\item $C_I(d)$ is the constant appearing in the standard regularity estimate (\ref{vharmonic}).
\item $C_{Mean}(d,\lambda)$ denotes the constant from the mean-value property in Theorem 2, i.e. (\ref{mvp}).
\end{itemize} 
All of these constants are assumed to be greater than 1.\\

Recall that for $i \in \left\{1,..., d-1\right\}$ we would like to construct the half-space-adapted corrector $\hp_{b_i}$ satisfying (\ref{halfspacecorrector}) and the corresponding skew-symmetric vector potential. We do this by ``correcting'' the whole-space corrector/vector potential pair that is assumed to exist in Theorem \ref{existenceofcorrectors} and showing that the corrections are sublinear thanks to (\ref{reqwsc}). In the following steps we construct the corrections, which must be built iteratively on increasingly large dyadic annuli. In particular, the main idea is to build a sublinear half-space-adapted corrector/vector potential pair up to a certain scale (Steps 1 and 2), use this pair along with Theorem \ref{excessdecay} to obtain a regularity theory up to this scale, apply the regularity theory to build a sublinear half-space-adapted corrector/ vector potential pair on a larger scale (Step 3), and then pass to the limit in this process (Step 4).

First, we fix an arbitrary initial radius $r_0\geq1$ and introduce two sets of functions:
\begin{enumerate}[topsep=0pt]
 \item[i)] Let $\{\eta_n\}_{n\geq-1}$ be a smooth radial partition of unity subordinate to the covering of $\mathbb{R}^d$ by $\{ B_{r_02^{n+1}} \setminus B_{r_0 2^{n-1}} \}_{n \geq 0} \cup B_{r_0}$ such that $|\nabla \eta_n | \leq \frac{4}{r_0 2^{n}}$.
 \item[ii)]  For each set in the cover we define a smooth one-dimensional cut-off function $L_n(x) = L_n(x_d)$ satisfying $|L_n(x_d)| = 1$ for  $|x_d| \leq l_n$  and $L_n(x_d) = 0$ for $|x_d| \geq 2 l_n$ such that $|\nabla L_n| \leq \frac{2}{l_n}$. 
\end{enumerate}
Specific values for the heights $l_n$ are chosen in the proof of Lemma \ref{energyvarphi}. For any radius $r>0$, to measure the size of the corrections to $\phi_{b_i}$ and $\sigma_{b_i}$ on $B_r^+$, we split the dyadic annuli into two groups: near-field, when $16r >r_02^{n+1}$,  and far-field, when $16r \leq r_02^{n+1}$.\\

\noindent \textbf{Step 1: Estimate for the near-field contributions of the correction $\varphi_{b_i}$ to the whole-space corrector $\phi_{b_i}$ when $i \neq d$.}\\

The correction to $\phi_{b_i}$, which we will call $\varphi_{b_i}$, that will enforce the desired boundary condition is a weak solution of
\begin{subequations}
 \label{correctionphi}
 \begin{align}
  ~~~~ -\di (a \nabla \varphi_{b_i}) & =0 && \textrm{in } \quad \hs,\\
  \label{correctionphibc}
  ~~~~ e_d \cdot a \nabla \varphi_{b_i} & = - e_d  \cdot a \nabla ( \phi_{b_i} + b_i \cdot x) && \textrm{on } \quad \bhs
 \end{align}
\end{subequations}
where the class of test functions is given by $H^1_{bdd}(\hs)$. For the boundary condition \eqref{correctionphibc} recall from \eqref{wsceqp} that the current $a \nabla (\phi_{b_i} + b_i \cdot x)$ is a solenoidal field, which means that $e_d \cdot  a \nabla (\phi_{b_i} + b_i \cdot x)$ has a trace in $H^{-1/2}(\bhs)$ and we may interpret the boundary condition in the distributional sense. 

To solve \eqref{correctionphi} and enforce the sublinearity of the correction, we split $\hs$ into the dyadic annuli introduced above (indexed by $n$) and for each $n \in \{ -1,0,1,..\}$ seek a solution $\varphi_{b_i}^n$ to 
\begin{subequations}
 \label{correctionphi2}
 \begin{align}
  ~~~~ -\di (a \nabla \varphi_{b_i}^n) & =0 && \textrm{in } \quad \hs,\\
  \label{correctionphibc1}
  ~~~~ e_d \cdot a \nabla \varphi_{b_i}^n & = - \eta_n e_d \cdot  a \nabla ( \phi_{b_i} + b_i \cdot x) && \textrm{on } \quad \bhs.
 \end{align}
\end{subequations}
The ansatz for the correction is then $ \varphi_{b_i} = \sum_{n=-1}^{\infty} \varphi_{b_i}^n$, which makes the ansatz for the half-space-adapted corrector $\hp_{b_i} =  \phi_{b_i} + \sum_{n=-1}^{\infty} \varphi_{b_i}^n $.
 
For a fixed $n\geq -1$ a solution to \eqref{correctionphi2} can be found in $\dot{H}^1(\hs)$. In particular, we perform a Lax-Milgram argument in the space $L^{2d/(d-2)}(\hs)\cap \dot{H}^1 (\hs)$ endowed with the inner-product $\langle \cdot, \cdot \rangle_{\dot{H}^1}$ to find a function $\tilde{\varphi}_{b_i}^n$ that satisfies the weak formulation
\begin{align}
 \label{weakformulationcorrect1}
 \int_{\hs} \nabla u \cdot a \nabla \tilde{\varphi}_{b_i}^n \, dx = - \int_{\bhs} u \eta_n e_d \cdot a \nabla(\phi_{b_i} + b_i \cdot x) \, dS
\end{align}
for all test functions $u$ from the Lax-Milgram space. Notice first of all that the inclusion $H^1_{bbd}(\hs) \subseteq  L^{2d/(d-2)}(\hs)\cap \dot{H}^1 (\hs)$ is due to the critical Sobolev embedding and that \eqref{weakformulationcorrect1} is the weak formulation of \eqref{correctionphi2} for $u \in H^1_{bdd}(\hs)$. We then let $\varphi_{b_i}^n = \tilde{\varphi}_{b_i}^n - \fint_{B_1^+} \tilde{\varphi}_{b_i}^n \, \textrm{d}x$. 

For the actual Lax-Milgram argument we notice:
\begin{enumerate}[topsep=0pt]
 \item [i)] The space $L^{2d/(d-2)}(\hs)\cap \dot{H}^1 (\hs)$ endowed with $\langle \cdot, \cdot \rangle_{\dot{H}^1}$ is complete and, therefore, a Hilbert space thanks to the Sobolev embedding.
 \item [ii)]The integral on the right-hand side of \eqref{weakformulationcorrect1} is well-defined due to the compact support of $\eta_n$.
 \item [iii)]The bilinear form on the left-hand side of \eqref{weakformulationcorrect1} is coercive due to the uniform ellipticity of $a$.
 \item [iv)] We then check that the right-hand side defines a bounded operator on the Lax-Milgram space. In particular, we notice that
\begin{align*}
 & \quad -\int_{\bhs} u \eta_n e_d \cdot a \nabla(\phi_{b_i} + b_i \cdot x) \, dS\\
  \quad \quad \quad \quad = ~~~~&\int_{\hs} \nabla (u \eta_n)  \cdot a \nabla(\phi_{b_i} + b_i \cdot x) \, dx\\
 \quad \leq ~~~~& C(d, \lambda, n) \left( \int_{\hs} | \nabla u |^2 \, dx \right)^{1/2} \left( \int_{B_{r_02^{n+1}}} | \nabla (\phi_{b_i} + b_i \cdot x) |^2 \, dx \right)^{1/2},
\end{align*}
where we have used the compact support of $\eta_n$, the boundedness of $a$, the critical Sobolev embedding, and that $\phi_{b_i} + b_i \cdot x$ is $a$-harmonic. As $\phi_{b_i} + b_i \cdot x \in H^1_{loc}(\mathbb{R}^d)$ it follows that the right-hand side of \eqref{weakformulationcorrect1} is a bounded operator.
\end{enumerate}
Having checked all of the criterion, we find that we may apply Lax-Milgram to obtain a solution $\tilde{\varphi}_{b_i}^n \in L^{2d/(d-2)}(\hs)\cap \dot{H}^1 (\hs)$ to \eqref{weakformulationcorrect1}.

It then remains to show that the ansatz $\sum_{n=-1}^{\infty} \varphi_{b_i}^n$ converges and is sublinear. For this purpose we introduce the notation $\hp_{b_i, N} =    \phi_{b_i}+\sum_{n=-1}^N \varphi_{b_i}^n$ and notice that the solution $\varphi_{b_i}^n$  of \eqref{correctionphi2} also solves
\begin{subequations}
    \label{etavarphi}
    \begin{align}
    \label{etavarphiequation}
    ~~~~ -\di (a ( \nabla \varphi_{b_i}^n + \eta_n L_n \nabla ( \phi_{b_i} + b_i \cdot x)) &= - \di (\eta_n L_n a \nabla ( \phi_{b_i} + b_i \cdot x)) &&\textrm{in }\quad \hs.\\
    \label{etavarphiboundary}
    ~~~~ e_d \cdot a ( \nabla \varphi_{b_i}^n + \eta_n L_n \nabla ( \phi_{b_i} + b_i \cdot x)) & =  0 &&\textrm{on }\quad \bhs,
    \end{align}
    \end{subequations} 
where we have smuggled in the vertical cut-off $L_n$.

The main observation of Step 1 is that we may choose the heights of the supports of the vertical cut-off functions $L_n$ such that the energy estimate for (\ref{etavarphi}) provides a sufficient bound for the size of the near-field contributions. In particular, we see in Step 3 that this allows us to build a sublinear half-space-adapted corrector/vector potential pair on $B_{8r_0}^+$.

\begin{customlemma}{2.1}(Energy estimate for $\varphi_{b_i}^n$)\label{energyvarphi} Assume that the conditions in Theorem 1 are satisfied. Then there exists $C_{1}(d, \lambda) \geq 1$ such that for each $n\geq -1$ there is a height $l_n>0$ so that for any $r >0$ and $i \in \{1,..,d-1\}$ it holds that:
\begin{align}\label{energyfirst}
 \left( \fint_{B_r^+} |\nabla \varphi_{b_i}^n |^2 \, dx \right)^{1/2} \leq C_1(d, \lambda) \left(\frac{r_02^{n+1}}{r}\right)^{d/2} \delta_{r_02^{n+1}}^{1/3}.
\end{align}
\noindent  In particular, when $ 16 r > r_0 2^{n+1}$ we have that 
\begin{align}
 \label{assumption1}
\left( \fint_{B_r^+}|\nabla \varphi_{b_i}^n|^2 \,dx  \right)^{1/2} \leq C_2(d, \lambda) \min \bigg\{1, \left(\frac{r_0 2^{n+1}}{r}\right)^{d/2} \bigg\} \delta_{r_02^{n+1}}^{1/3}
\end{align}
\noindent with $C_2 := C_{Mean} C_1 8^{d}$.
\end{customlemma}

\begin{proof}[Proof of Lemma \ref{energyvarphi}] 

For convenience, in this proof we use the notation $R = r_02^{n+1}$. As the desired inequalities \eqref{energyfirst} and \eqref{assumption1} for $\varphi_{b_i}^n$ only involve $\nabla \varphi^n_{b_i}$ and \eqref{etavarphi} does not see the subtraction of constants from $\varphi_{b_i}^n$, we assume for this argument that $\fint_{B_R^+} \varphi_{b_i}^n \, dx = 0$.  Testing \eqref{etavarphi} with $\varphi_{b_i}^n$ yields
\begin{align}
\label{energyvarphi1}
\begin{split}
\int_{\hs} \nabla \varphi_{b_i}^n  \cdot a \nabla \varphi_{b_i}^n \,dx = &  - \int_{\hs} \eta_n L_n \nabla \varphi_{b_i}^n \cdot a \nabla (\phi_{b_i} + b_i \cdot x)\,dx\\
&- \int_{\hs} \varphi_{b_i}^n L_n \nabla \eta_n  \cdot a \nabla (\phi_{b_i} + b_i \cdot x)\,dx \\
& - \int_{\hs} \varphi_{b_i}^n \eta_n \nabla L_n \cdot a \nabla (\phi_{b_i} + b_i \cdot x)\,dx.
\end{split}
\end{align}
We first treat the third term on the right-hand side of (\ref{energyvarphi1}). As $ i \in \left\{1,..., d-1 \right\}$ we have that $e_d \cdot a_{hom} b_i =0$, which, thanks to (\ref{wsceqs}), enables us to express  $e_d \cdot a \nabla (\phi_{b_i} + b_i \cdot x)$ in divergence form. In particular, we know that the equation
\begin{align}
\nabla_k \cdot \sigma_{b_i d k} =  e_d \cdot a \nabla (\phi_{b_i} + b_i \cdot x) \quad \quad \quad \textrm{ in } \quad \mathbb{R}^d
\end{align}
is satisfied in the distributional sense. Due to the choice of $L_n$ and $\eta_n$ we have that $\varphi_{b_i}^n \eta_n \nabla L_n \in H^1_0(\hs)$.

In the following computation $k$ indexes the entries of the vector $\sigma_{b_i d}$. Making use of the identity  $\nabla L_n =   \partial_d L_n e_d$ we write:
\begin{align}
\label{energyvarphi2}
\begin{split}
 \int_{\hs} \varphi_{b_i}^n \eta_n \nabla L_n \cdot a \nabla (\phi_{b_i} + b_i \cdot x) \,dx=  & \int_{\hs} \varphi_{b_i}^n \eta_n \partial_d L_n e_d \cdot a \nabla (\phi_{b_i} + b_i \cdot x)\,dx\\
  =  & - \int_{\hs}  \partial_d L_n  \partial_k (\varphi_{b_i}^n \eta_n)  \sigma_{b_i d k}\,dx\\
     & - \int_{\hs} \varphi_{b_i}^n \eta_n \partial_d^2 L_n  \sigma_{b_i d d}\,dx\\
  =  & - \int_{\hs}  \partial_d L_n  \partial_k (\varphi_{b_i}^n \eta_n) \sigma_{b_i d k}\,dx.
 \end{split}
\end{align}
\noindent Notice that the boundary terms in the integration by parts vanish and the last equality follows from the skew-symmetry of $\sigma_{b_i}$.

Making use of (\ref{energyvarphi2}), the uniform ellipticity and boundedness of $a$, and the Poincar\'{e} inequality with zero-average on $B_{R}^+$, we find that (\ref{energyvarphi1}) implies
\begin{align}
\label{energyvarphi3}
\begin{split}
\int_{\hs} |\nabla  \varphi_{b_i}^n|^2\,dx  \lesssim & \left( \int_{B_R^+} |\nabla \varphi_{b_i}^n|^2 \,dx \right)^{1/2} \left( \left( \int_{\textrm{supp}(\eta_n L_n)}  |\nabla (\phi_{b_i} + b_i \cdot x)|^2 \,dx \right)^{1/2} \right.\\
& + \left( \int_{B_R^+}  |\partial_d L_n  \eta_n  \sigma_{b_i  d}|^2 \,dx \right)^{1/2} \\
& + R  \left( \int_{\textrm{supp}(\eta_n L_n)} |\nabla \eta_n|^2 |\nabla (\phi_{b_i} + b_i \cdot x)|^2 \,dx \right)^{1/2}\\
& + \left.  R \left( \int_{B_R^+}  |\partial_d L_n  \partial_k \eta_n  \sigma_{b_i d k}|^2 \, dx \right)^{1/2} \vphantom{\left( \int_{\textrm{supp}(\eta_n L_n)}  | \nabla (\phi_{b_i} + b_i \cdot x)|^2 \right)^{1/2}} \right).\\
\end{split}
\end{align}
\noindent We can simplify this expression by recalling that $|\nabla \eta_n| \leq \frac{8}{R}$ and $ | \partial_d L_n | \leq \frac{2}{l_n}$, which allows us to write
\begin{align}
 \label{energyvarphi4}
 \begin{split}
\int_{\hs} |\nabla  \varphi_{b_i}^n|^2 \,dx  \lesssim & \left( \int_{B_R^+} |\nabla \varphi_{b_i}^n|^2 \,dx \right)^{1/2} \left( \left( \int_{\textrm{supp}(\eta_n L_n)}  |\nabla (\phi_{b_i} + b_i \cdot x)|^2 \,dx \right)^{1/2} \right.\\
& + \left.\frac{1}{l_n} \left( \int_{B_R^+}  |\sigma_{b_i d}|^2 \,dx \right)^{1/2} \vphantom{\left( \int_{\textrm{supp}(\eta_n L_n)}  | \nabla (\phi_{b_i} + b_i \cdot x)|^2  \right)^{1/2}} \right).
 \end{split}
\end{align}

We must still treat the first term on the right-hand side of (\ref{energyvarphi4}), which we do with a box-wise Caccioppoli estimate. In particular, we cover $\textrm{supp}(\eta_n L_n)$ with cubes of side length $4l_n$. If we denote the $d$-dimensional cube with center $z\in \mathbb{R}^d$ and side length $l\in \mathbb{R}$ by $C_l(z)$, we may find a set of points
\begin{align}
\label{defnofs}
\begin{split}
 S  = & \left\{z \in \mathbb{R}^d \hspace{.2cm} \middle| \hspace{.2cm} |\textrm{supp}(\eta_n L_n) \setminus \cup_{z \in S }C_{4l_n}(z) | = 0 \right.\\
 & \left. \textrm{ and for all } x \in \mathbb{R}^d \textrm{ we have that } \sum_{z\in S} \chi_{C_{6l_n}(z)}(x) \leq 2^d \vphantom{z \in \mathbb{R}^d} \right\}.
\end{split}
\end{align}
Then, for each box $C_{4l_n}(z)$ we let $\tilde{C}_{4 l_n, 6l_n,z}$ denote the smooth cut-off of $C_{4l_n}(z)$ in the box of side length $6l_n$ centered around it. In particular, we require that
\begin{align}
\label{defnCjcutoff}
\tilde{C}_{4l_n, 6l_n,  z } (x) = \bigg\{
\begin{array}{ll}
       1 & \textrm{  if  } x \in C_{4l_n}(z) \\
       0 & \textrm{  if  } x \notin C_{6l_n}(z)\\
\end{array}
\end{align}
\noindent and that $|\nabla \tilde{C}_{4 l_n, 6l_n, z}| \leq 1/l_n$.

For each $z \in S$ we test equation (\ref{wsceqp}) with $(\tilde{C}_{4l_n,6l_n, z})^2 \eta_n^2 (\phi_{b_i} + b_i \cdot (x-z))$. After using the uniform ellipticity and boundedness of $a$ and Young's inequality we obtain that
\begin{align}
\label{varphienergycac2}
 \int_{C_{4l_n}(z)\cap B_R} |\nabla(\phi_{b_i} + b_i \cdot x)|^2 \,dx \lesssim \frac{1}{l_n^2} \int_{C_{6l_n}(z) \cap B_R}  |\phi_{b_i} + b_i \cdot (x - z) |^2\,dx,
\end{align}
where we have also used that we may choose $l_n$ to satisfy $l_n \ll R$. Breaking up $\textrm{supp}(\eta_n L_n)$ into the cubes $C_{4l_n}(z)$ with centers $z \in S$ and applying (\ref{varphienergycac2}) on each cube, we find that
\begin{align}
\label{varphienergycac3}
 \begin{split}
  \int_{\textrm{supp}(\eta_n L_n)}  |\nabla (\phi_{b_i} + b_i \cdot x)|^2 \,dx & \leq \sum_{z\in S} \int_{C_{4l_n}(z)\cap B_R} |\nabla(\phi_{b_i} + b_i \cdot x)|^2 \,dx\\
  &\lesssim  \frac{1}{l_n^2} \displaystyle\sum_{z \in S} \int_{C_{6l_n}(z)\cap B_R} |\phi_{b_i} + b_i \cdot (x - z) |^2\,dx\\
  &\lesssim \frac{1}{l_n^2}\left(  \int_{B_R} |\phi_{b_i}|^2 \,dx + \sum_{z \in S} \int_{C_{6l_n}(z)\cap B_R} |x - z|^2\,dx\right)\\
  &\lesssim R^d \left(\frac{R}{l_n}\right)^2 \delta_{R}^2  + R^{d-1}l_n.
 \end{split}
\end{align}
\noindent Here we have used that the longest diagonal in a $d$-dimensional box of side length $6l_n$ has length $6 l_n d^{1/2}$ and also \eqref{different_basis}.

Combining (\ref{varphienergycac3}) with (\ref{energyvarphi4}) we find that for all $r>0$ it holds that
\begin{align}
\fint_{B_r^+} |\nabla  \varphi_{b_i}^n|^2 \,dx& \lesssim \left(\frac{R}{r}\right)^d \left( \left(\frac{R}{l_n}\right)^2 \delta_{R}^2  + \frac{l_n}{R} \right).
\end{align}
\noindent Letting $l_n = \alpha R$ and plugging in the optimal $\alpha = \delta_R^{2/3}$ yields (\ref{energyfirst}). Lastly, we note that (\ref{assumption1}) is a trivial consequence of (\ref{energyfirst}) and that $C_{Mean} \geq 1$.
\end{proof}

\noindent \textbf{Step 2: Estimate for the correction $\psi_{b_i}$ of the whole-space vector potential $\sigma_{b_i}$ when $i \neq d$.}\\

In this step we construct intermediate half-space-adapted vector potentials, which we call $\hsi_{b_i, N}$, corresponding to the $\hp_{b_i, N}$ from the last step. In particular, we construct each $\hsi_{b_ij, N}$ to be a distributional solution of
\begin{align}
    \label{partialhalfspacepotential}
    \nabla_k \cdot \hsi_{b_ijk, N} = e_j \cdot ( a(b_i + \nabla \hp_{b_i, N} ) - a_{hom}b_i) \quad \quad \quad\textrm{in} \quad \hs.
\end{align}
Our strategy for this construction is to correct the whole-space $\sigma_{b_ijk}$  with a modification $\psi_{b_i jk, N}$ that satisfies
\begin{align}
\label{sigmacorrection}
\nabla_k \cdot \psi_{b_ijk,N} = e_j \cdot (a(b_i + \nabla \hp_{b_i, N}) - a(b_i + \nabla \phi_{b_i})) \quad \quad \quad \textrm{in} \quad \hs.
\end{align} 
Taking the ansatz $\hsi_{b_ijk,N} = \sigma_{b_ijk } + \psi_{b_ijk, N}$, we must then ensure that $\psi_{b_ijk,N}$ is sublinear and skew-symmetric in $j$ and $k$. 

To obtain the desired corrections $\psi_{b_i jk, N}$ we again decompose $\mathbb{R}^d$ into the dyadic annuli from the last step and then consider the solutions $v_{b_ij, N}^n : \hs \rightarrow \mathbb{R}$ of
\begin{subequations}
\label{equnvn}
\begin{align}
\label{equnvna}
-\Delta v_{b_i j,N}^n &= \di (\eta_n x_j a (\nabla \hp_{b_i, N} - \nabla \phi_{b_i}) ) \quad &&\textrm{in } \hs,
\\
\label{vnbc1}
v_{b_i j,N}^n &= 0 \quad &&\textrm{for $j \neq d$  on } \bhs,
\\
\label{vnbc2}
\partial_d v_{b_i d,N}^n &= 0 \quad &&\textrm{on } \bhs.
\end{align}
\end{subequations}
In particular, for fixed $n$ and $j \neq d$ we will find a solution $v_{b_i j,N}^n \in \dot{H}^1(\hs)$, which satisfies \eqref{equnvna} in the distributional sense and satisfies the boundary condition \eqref{vnbc1} in the trace sense; when $j =d$ we find a solution $v_{b_i j, N}^n \in \dot{H}^1(\hs)$ which satisfies \eqref{equnvna} and \eqref{vnbc2} in the weak sense with the class of test function being given by $H^1_{bdd}(\hs)$. We plan to find these solutions using a Lax-Milgram argument.

We first consider the case $j \neq d$. Using a Lax-Milgram argument in the space $\dot{H}^1_0(\hs)$ we are able to find $v_{b_i j, N}^n \in \dot{H}^1_0(\hs)$ satisfying the weak formulation 
\begin{align}
\label{weak_form_jneqd}
 \int_{\hs} \nabla u \cdot \nabla v_{b_i j,N}^n \textrm{d} x = - \int_{\hs} \eta_n x_j   \nabla u \cdot a (\nabla \hp_{b_i, N} - \nabla \phi_{b_i})\textrm{d} x
\end{align}
for any $u \in \dot{H}^1_0(\hs)$. For this, notice that the bilinear form on the left-hand side of \eqref{weak_form_jneqd} is coercive and that, as $\nabla \phi_{b_i,N}^{\h} - \nabla \phi_{b_i} = \sum_{n=-1}^N \nabla \varphi_{b_i}^n$, it follows from $\textrm{supp}(\eta_n) \subset B_{r_02^{n+1}}^+$ and Lemma \ref{energyvarphi} that $\eta_n x_j a(\nabla \phi_{b_i,N}^{\h} - \nabla \phi_{b_i}) \in L^2(\hs)$;  So, the right-hand side of \eqref{weak_form_jneqd} defines a bounded operator on $\dot{H}^1_0(\hs)$. The equality \eqref{weak_form_jneqd} is the weak formulation of \eqref{equnvn} for $j \neq d$ when $u \in C^{\infty}_c(\hs)$.

For the case $j=d$ we use the space $L^{2d/(d-2)}(\hs)\cap \dot{H}^1 (\hs)$ endowed with the inner-product $\langle \cdot, \cdot \rangle_{\dot{H}^1}$ for our Lax-Milgram argument. It is clear from the discussion for the case $j \neq d$ above that we may find $v_{b_i d,N}^n \in L^{2d/(d-2)}(\hs)\cap \dot{H}^1 (\hs)$ such that the weak formulation \eqref{weak_form_jneqd} is satisfied for all $u \in L^{2d/(d-2)}(\hs)\cap \dot{H}^1 (\hs)$. We additionally notice that for $j = d$ the equality \eqref{weak_form_jneqd} is the weak formulation of \eqref{equnvn} for $u \in H^1_{bdd}(\hs)$. In particular, when we test $\eqref{equnvn}$ for $j =d$ with $u\in H^1_{bdd}(\hs)$ the boundary term on $\bhs$ vanishes due to the boundary condition \eqref{vnbc2} and since $x_d = 0$; and the boundary term at infinity vanishes since the test function $u$ has compact support.

We would then, for each $i \in \left\{1,..., d-1 \right\}$, $j \in \left\{1,...,d \right\}$ and $N \geq -1$, like to sum over all of the contributions $v_{b_ij,N}^n$. In order to ensure that this sum converges we subtract-off the initial linear growth of each $v_{b_ij,N}^n$. For this purpose, set
\begin{align}
\label{DefinitionbjMn}
c_{b_i j,N}^n = :\bigg\{
\begin{array}{ll}
       0 & \textrm{  if  } n=-1 \\
     \nabla v_{b_i j,N}^n(0) & \textrm{  if  } n \neq -1
\end{array}
\end{align}
and, furthermore, let
\begin{align}
\label{Definitiond}
d_{b_i j,N}^n = :\bigg\{
\begin{array}{ll}
       0 & \textrm{  if  } j\neq d \\
     \fint_{B_1^+} (v_{b_i d, N}^n - c_{b_i d, N}^n \cdot x)\, \textrm{d}x & \textrm{  if  } j =d.
\end{array}
\end{align}
It is shown in Lemma \ref{2.2} that $v_{b_i j, N}^n - c_{b_i j,N}^n \cdot x - d_{b_i j,N}^n$ has quadratic behavior inside $B_r^+$ whenever $16r \leq r_02^{n+1}$ (i.e. when $\varphi_n$ is a far-field contribution). This is then used in Lemma \ref{psiestimate} to establish that for any $r>0$ the expression $\sum_{n=-1}^{\infty} ( v_{b_i j, N}^n - c_{b_i j,N}^n \cdot x- d_{b_i j,N}^n)$ converges absolutely in $H^1(B_r^+)$ for all $j \in \left\{1,..., d\right\}$. We denote $v_{b_i j, N} = \sum_{n=-1}^{\infty}  (v_{b_i j, N}^n - c_{b_i j,N}^n \cdot x- d_{b_i j,N}^n)$.

We observe that the sum $v_{b_i, N}$ satisfies
\begin{subequations}
\label{equnv}
\begin{align}
\label{equnva}
-\Delta v_{b_i j,N} &= \di (x_j a (\nabla \hp_{b_i, N} - \nabla \phi_{b_i}) ) \quad &&\textrm{in } \hs,
\\
\label{bccorrection1}
v_{b_i j,N} &= 0 \quad &&\textrm{for $j \neq d$  on } \bhs,
\\
\label{bccorrection2}
\partial_d v_{b_i d,N} &= 0 \quad &&\textrm{on } \bhs.
\end{align}
\end{subequations}
In particular, using the definitions \eqref{DefinitionbjMn} and \eqref{Definitiond}, we find that for every $n \geq -1$ if $j \neq d$ then $v_{b_ij,N}^n - c_{b_ij,N}^n \cdot x- d_{b_i j,N}^n$ is a distributional solution of \eqref{equnvn} and $v_{b_id,N}^n - c_{b_id,N}^n \cdot x- d_{b_i d,N}^n$ is a weak solution with test functions taken from $H^1_{bdd}(\hs)$.

Differentiating (\ref{equnv}) we find that $\nabla_k \cdot v_{b_i k, N}$ is harmonic on $\hs$ with homogeneous Dirichlet boundary conditions. Assuming for now that $\nabla_k \cdot v_{b_i k, N}$ is sublinear (see Lemma 2.3), the zeroth-order Liouville principle for harmonic functions with homogeneous Dirichlet boundary conditions then implies that $\nabla_k \cdot v_{b_i k, N} \equiv 0$ on $\hs$. This, in combination with \eqref{equnva}, allows us to conclude that the ansatz
\begin{align}
\label{ansatzpsi}
 \psi_{b_ijk, N} = \partial_k v_{b_i j, N}-  \partial_j v_{b_i k,N}
\end{align}
solves (\ref{sigmacorrection}). For the complete details please see the calculation in Step 2 of Section 2.2 of \cite{FischerRaithel}. To conclude, we notice that \eqref{ansatzpsi} is antisymmetric in $j$ and $k$ and prove that it is sublinear. 

\begin{customlemma}{2.2}\label{2.2}  Assume that the conditions in Theorem 1 are satisfied and the heights $l_n$ are chosen as specified in Lemma 2.1. Let $N \geq -1$ and $n\geq -1$. Then for $j,k \in \left\{1,...,d \right\}$ and  $i \in \left\{1,...,d-1 \right\}$ it holds that 
\begin{align}
\label{assumption1.2}
\begin{split}
&\frac{1}{r}\left( \fint_{B_r^+}|\partial_k (v_{b_i j,N}^n - c_{b_i j,N}^n \cdot x )|^2 \,dx \right)^{1/2}\\
& \leq C_3(d) \min\left\{1, \frac{r_02^{n+1}}{r}\right\}\left(\fint_{B_{r_02^{n+1}}^+}|\nabla \hp_{b_i,N} - \nabla \phi_{b_i}|^2 \,dx\right)^{1/2}
\end{split}
\end{align}
\noindent for all $r >0$ and $C_3(d) := 16 C_4(d) C_I(d)$, where the constant $C_4(d)$ is specified in the proof. 
\end{customlemma}

\begin{proof} Our argument mainly relies on two observations:

\begin{itemize}
\item[i)] Testing the weak formulation \eqref{weak_form_jneqd} with $v_{b_ij,N}^n$ and combining the resulting energy estimate with
\begin{align}
|c_{b_i j,N}^n|\leq C(d) \left( \fint_{B_{R/4}^+} |\nabla v_{b_i j,N}^n|^2 \, \textrm{d}x  \right)^{1/2},
\end{align}
which is a result of the mean-value property for harmonic functions and $\partial_k v_{b_i j,N}^n$ being harmonic with homogeneous Dirichlet or Neumann boundary data (depending on $j$ and $k$), gives that for any two radii $R>r>0$ such that $r \geq \frac{1}{16}R$ the relation
\begin{align}
\label{energyvn}
\left( \fint_{B_r^+}|\nabla v_{b_i j,N}^n - c_{b_i j,N}^n|^2 \,dx\right)^{1/2} & \leq C_4(d) R \left( \fint_{B_R^+}|\nabla \hp_{b_i,N} - \nabla \phi_{b_i}|^2  \,dx \right)^{1/2}
\end{align}
holds for some $C_4(d)\geq 1$.
\item[ii)] For any $R>0$, if $w \in L^2(B_R^+)$ is harmonic with either homogeneous Neumann or homogeneous Dirichlet boundary conditions on $\fpb_{R}$ then for $ r \in (0, R/4]$ we have
\begin{align}
\label{vharmonic}
\left(\fint_{B_r^+}|w - w(0)|^2 \,dx \right)\leq C_I(d)\frac{r}{R}\left(\fint_{B_{R}^+} |w|^2 \,dx\right)^{1/2}.
\end{align}
This follows, e.g., from an iterative application of an appropriate Caccioppoli estimate.\\

\noindent \textbf{Remark}: Here we use that, by iterating the usual Caccioppoli estimate for the Dirichlet case and that from Lemma \ref{Caccioppolifornoflux} for the Neumann case, one finds that $w \in H^{\lfloor d/2 \rfloor +2}(B_{R/2}^+)$. As $H^{\lfloor d/2 \rfloor +2}(B_{R/2}^+)$ embeds into $C^{1,\gamma}\left(\overline{B}_{R/2}^+ \right)$ for some $1>\gamma>0$, this, in particular, means that it makes sense to write $w(0)$.\\
\end{itemize} 

To obtain Lemma 2.2 for a fixed radius $r>0$ one considers near-field and far-field contributions separately. The desired estimate \eqref{assumption1.2} follows immediately from \eqref{energyvn} for the near-field contributions, i.e. when $ r \geq r_02^{n-3}$. To treat the far-field contributions we notice that $\partial_k v_{b_i j,N}^n - e_k \cdot c_{b_i j,N}^n$  is harmonic on $B_{r_02^{n-1}}^+$ with homogeneous Dirichlet or homogeneous Neumann boundary conditions on $\fpb_{r_02^{n-1}}$ (depending on $j$ and $k$). As $ r \leq r_02^{n-3}$ we may first apply (\ref{vharmonic}) and then (\ref{energyvn}) in the following way:
\begin{align}
\label{lemma22farfield}
\begin{split}
&\frac{1}{r}\left(\fint_{B_r^+}|\partial_k (v_{b_i j,N}^n - c_{b_i j,N}^n \cdot x)|^2 \, dx  \right)^{1/2}\\
& \leq C_I\frac{1}{r_02^{n-1}}\left( \fint_{B_{r_02^{n-1}}^+}|\partial_k (v_{b_i j,N}^n - c_{b_ij,N}^n\cdot x)|^2\,dx\right)^{1/2}\\
& \leq 4 C_4 C_I \left(\fint_{B_{r_02^{n+1}}^+}|\nabla \hp_{b_i,N} - \nabla \phi_{b_i}|^2 \,dx \right)^{1/2}.
\end{split}
\end{align}
\end{proof}

\begin{customlemma}{2.3} \label{psiestimate} Assume that the conditions in Theorem 1 are satisfied and the heights $l_n$ are chosen as specified in Lemma 2.1. Let $i \in \{1,..., d-1\}$, $j \in\left\{1,...,d \right\}$, and $N \geq -1$. Then:
\begin{enumerate}
 \item[i)] The sum $\displaystyle\sum_{n=-1}^{\infty} ( v_{b_i j,N}^n - c_{b_i j,N }^n \cdot x - d_{b_i j,N}^n)$ converges absolutely in $H^1(B_r^+)$ for any $r>0$. We denote the limit by $v_{b_i j,N}$.
 \item[ii)] The expression $\nabla_k \cdot v_{b_i k, N}$ is sublinear in the sense that 
 \begin{align*}
 \displaystyle\lim_{r\rightarrow \infty}\frac{1}{r}\left( \fint_{B_r^+}|\nabla_k \cdot v_{b_i k, N}|^2 \, \textrm{d}x \right)^{1/2} = 0.
\end{align*}
 \item[iii)] For $r >0$ and $k\in \{1,...,d\}$ the ansatz $\psi_{b_i jk,N}$ satisfies the estimate 
\begin{align}
\label{estimatepsi}
\begin{split}
\frac{1}{r}\left( \fint_{B_r^+}| \psi_{b_i jk,N}|^2 \,dx \right)^{1/2} & \leq 2 C_3(d, \lambda) \sum_{n=-1}^{\infty}\min \left\{1, \frac{r_02^{n+1}}{r} \right\}\\
& ~~~~~~~~~ \times \left( \fint_{B_{r_02^{n+1}}^+} |\nabla \hp_{b_i,N} - \nabla \phi_{b_i}|^2 \,dx \right)^{1/2}.
\end{split}
\end{align}
\end{enumerate}
\end{customlemma}

\vspace{.25cm}

\begin{proof}
By Lemma \ref{2.2} it is clear that if 
\begin{align}
\label{lemma231}
\displaystyle\sum_{n=-1}^{\infty} \left( \fint_{B_{r_02^{n+1}}}|\nabla \hp_{b_i,N} - \nabla \hp_{b_i}|^2 \,dx\right)^{1/2}<\infty
\end{align} 
then $\sum_{n=-1}^{\infty} \nabla (v_{b_i j,N}^n - c_{b_i j,N }^n \cdot x)$ converges absolutely in $L^2(B_r^+)$ for any $j \in \left\{1,...,d \right\}$ and any $r >0$. Conveniently, (\ref{lemma231}) follows easily from the identity $\hp_{b_i,N} - \phi_{b_i} = \sum_{n=-1}^N \varphi_{b_i}^n$  and \eqref{energyfirst} from Lemma \ref{energyvarphi}.

When $j \neq d$ the homogeneous Dirichlet boundary data of $ v_{b_i j,N}^n - c_{b_i j,N }^n \cdot x - d_{b_i j,N}^n$ on $\bhs$ allows us to upgrade this to the absolute convergence of $\sum_{n=-1}^{\infty} (v_{b_i j,N}^n - c_{b_i j,N }^n \cdot x - d_{b_i j,N}^n)$ in $H^1(B_r^+)$. To treat the case $j =d$ we use that
\begin{align}
\fint_{B_1^+}  v_{b_i d,N}^n - c_{b_i d,N }^n \cdot x - d_{b_i d,N}^n \, dx = 0
\end{align}
for all $n \geq -1$. For any $u \in H^1(B_r^+)$ such that $\fint_{B_1^+}u \, \textrm{d}x =0$ notice that a combination of Jensen's inequality and the Poincar\'{e} inequality with zero average on $B_r^+$ yields
\begin{align}
\begin{split}
\label{averageball1poincare}
\int_{B_r^+}|u |^2 \, dx& = \int_{B_r^+}\left|u  - \fint_{B_1^+} u  \, dx \right|^2 \, dx\\
 \lesssim &  \, \int_{B_r^+} \left| u  - \fint_{B_r^+} u \, dx\right|^2  \, dx + r^d \int_{B_1^+}\left| u  - \fint_{B_r^+} u \, dx \right|^2 \, dx\\
 \lesssim &  \, r^{d+2} \int_{B_r^+}|\nabla u|^2 \, dx.
\end{split}
\end{align}
This, in particular, implies that we also in the case $j=d$ obtain absolute convergence in $H^1(B_r^+)$.

To finish, notice that the sublinearity of $\nabla_k \cdot v_{b_i k,N}$ follows from Lemma \ref{2.2} and the bound (\ref{lemma231}) using the dominated convergence theorem for sums. The estimate (\ref{estimatepsi}) for $\psi_{b_ij,N}$  also follows from Lemma \ref{2.2}.
\end{proof}

\noindent \textbf{Step 3: Inductive construction of $(\hp, \hsi)$ on larger scales.}\\

Notice that in the previous two steps the initial radius $r_0 \geq 1$ was arbitrary. In the current step we assume that the conditions of Theorem 1 are satisfied and choose a specific $r_0$, which is large enough so that for all $N \geq -1$ the intermediate half-space-adapted corrector/vector potential pair $(\hp_{,N}, \hsi_{,N})$ is sublinear in the sense that it satisfies condition (\ref{smallness}) from Theorem \ref{excessdecay} for $\alpha = 1/2$ and $r \geq r_0$. Furthermore, for this choice of $r_0$ and for $i \neq d$ we find that $\varphi_{b_i}^n$ satisfies (\ref{assumption1}) for any $r \geq r_0$ and $n \geq -1$. 

\begin{customlemma}{2.4}\label{chooser0} Assume that the conditions in Theorem 1 are satisfied and the heights $l_n$ are chosen as specified in Lemma 2.1. Then there exists a dyadic radius $r_0 = 2^{n_0}$ that does not depend on $N$ such that the following statements hold:\\

If for all $n \in \{-1,...,N \}$ and $i \in \{1,...,d-1\}$ it holds that 
\begin{align}
\label{assumption12}
 \left( \fint_{B_r^+} |\nabla \varphi_{b_i}^n |^2 \,dx \right)^{1/2} \leq C_2(d, \lambda) \min\left\{1, \left(\frac{r_02^{n+1}}{r}\right)^{d/2} \right\} \delta_{r_02^{n+1}}^{1/3}
\end{align}
for $r \geq r_0$ then $(\hp_{,N}, \hsi_{,N})$-- where we  let $\hp_{b_d,N} = \phi_{b_d}|_{\hs}$ and $\hsi_{b_d,N} = \sigma_{b_d}|_{\hs}$ for all $N \geq -1$-- satisfies condition (\ref{smallness}) from Theorem \ref{excessdecay} for $r \geq r_0$ and $\alpha =1/2$, i.e.
\begin{align}
\label{smallnessint}
\deltah_r(\hp_N, \hsi_N) \leq \frac{1}{C_{1/2}(d, \lambda)} \textrm{ for } r \geq r_0.
\end{align}\\
Furthermore, then (\ref{assumption12}) holds for $\varphi^{N+1}_{b_i}$ for all $i \in \{1,...,d-1\}$ and $r \geq r_0$.

\end{customlemma}

\begin{proof} Let $r>0$ be an arbitrary radius. Young's inequality and the Poincar\'{e} inequality with zero average then yield
\begin{align}
\label{intsublin1}
\begin{split}
     & \frac{1}{r} \left( \fint_{B_r^+} \sum_{i=1}^{d-1} \left( \left|\phi_{b_i} + \sum_{n=-1}^N \varphi_{b_i}^n - \fint_{B_r^+} (\phi_{b_i}+\sum_{n=-1}^N \varphi_{b_i}^n) \, dx \right|^2   + |\sigma_{b_i} + \psi_{b_i,N}|^2 \right) \, \textrm{d} x  \right.\\
     & \quad \quad \quad \quad \quad \quad\quad \quad \quad \quad \quad \quad \quad \quad\quad \quad \quad \quad \quad \quad \quad \quad \left. + \fint_{B_r} |\phi_{b_d}|^2 + |\sigma_{b_d}|^2  \,dx \vphantom{\sum_{n=-1}^N} \right)^{1/2}\\
& \quad \quad \quad  \quad \quad  \quad \quad\leq   \frac{4}{r}\left(\fint_{B_r}|(\phi_{b_i}, \sigma_{b_i})|^2 \, dx \right)^{1/2} + \, \, \sum_{i=1}^{d-1} \left(\frac{2}{r}\left(  \fint_{B_r^+} |\psi_{b_i,N}|^2 \, dx  \right)^{1/2} \right.\\ 
& ~~~~~~~~~~~~~~~~~~~~~~~~~~~~~~~~~~~~~~~~~~~~~~~~~~~~~~~~\left. +   C_P \sum_{n=-1}^N \left( \fint_{B_r^+}|\nabla \varphi_{b_i}^n |^2 \, dx\right)^{1/2}\right),
\end{split}
\end{align}
where we have made use of the Einstein summation convention. We then estimate the three terms on the right-hand side of (\ref{intsublin1}) separately. Notice that:
\begin{itemize}
\item[i)] Using \eqref{different_basis} we treat the first term as
\begin{align}
\label{intsublin2}
\frac{4}{r}\left(\fint_{B_r}|(\phi_{b_i}, \sigma_{b_i})|^2 \, dx\right)^{1/2}  & \leq 4 \left( \frac{d(d+1)}{2}\right)^{1/2}  \delta_{r}.
\end{align}
\item[ii)] For the second term, for any $i \in \left\{1,..., d-1 \right\}$ an application of Lemma \ref{psiestimate} and assumption \eqref{assumption12} yield that 
\begin{align}
\label{intsublin3}
\begin{split}
\frac{1}{r}\left( \fint_{B_r^+} |\psi_{b_i,N}|^2 \,dx \right)^{1/2} & \leq 2 d^2 C_3  \sum_{m=-1}^{\infty} \sum_{n=-1}^N \min \left\{1, \frac{r_02^{m+1}}{r} \right\}\\
& ~~~~~~~~~~~~~~~~~~ \times \left( \fint_{B_{r_02^{m+1}}^+} |\nabla \varphi_{b_i}^n|^2 \,dx  \right)^{1/2}\\
& \leq 2 d^2 C_3 C_2 \sum_{n=-1}^N \sum_{m=-1}^{\infty} \min\left\{ 1, 2^{d (n-m)/2}\right\}\delta_{r_02^{n+1}}^{1/3}\\
& \leq 2 d^2 C_3 C_2 \sum_{n=n_0}^{N+n_0+1} \left(n-n_0+\frac{1}{1-2^{-d/2}}\right)\delta_{2^n}^{1/3}.
\end{split}
\end{align}
\item[iii)] To treat the third term, for any $i \in \left\{1,..., d-1 \right\}$ assumption \eqref{assumption12} gives that 
\begin{align}
\begin{split}
\label{intsublin4}
 \sum_{n=-1}^N \left( \fint_{B_r^+}|\nabla \varphi_{b_i}^n |^2 \,dx \right)^{1/2} & \leq C_2\sum_{n=-1}^N \min \left\{1, \left( \frac{r_02^{n+1}}{r}\right)^{d/2} \right\}\delta_{r_02^{n+1}}^{1/3}\\
& \leq C_2 \sum_{n=n_0}^{N+n_0+1} \delta_{2^n}^{1/3}.
\end{split}
\end{align}
\end{itemize}

\noindent Combining these three estimates with \eqref{intsublin1} gives that 
\begin{align}
\begin{split}
\label{intsublin5}
\deltah_r(\hp_{,N}, \hsi_{,N}) & \leq 4\left( \frac{d(d+1)}{2}\right)^{1/2}\delta_r \\
& ~~~~~~~~~~~~~ + 4d^2C_3C_2 C_P \sum_{n=n_0}^{\infty} \left(n -n_0+ 1 + \frac{1}{1-2^{-d/2}}\right)\delta_{2^n}^{1/3}.
\end{split}
\end{align}
By our assumption \eqref{reqwsc} on the whole-space corrector/vector potential pair we find that we can choose the initial radius $r_0=2^{n_0}$ large enough, in a manner independent of $N$, such that \eqref{smallnessint} holds.

We then show that $\varphi_{b_i}^{N+1}$ satisfies \eqref{assumption12} for all $r \geq r_0$. Notice that \eqref{assumption12} for $r\geq r_0$ such that $\varphi_{b_i}^{N+1}$ is a near-field contribution, i.e. when $\frac{r_02^{N+2}}{r}\leq 16$, has already been shown in Lemma \ref{energyvarphi}. We, therefore, restrict ourselves to the case when $r \leq r_02^{N-2}$. By the argument above, the intermediate pair $(\hp_{,N}, \hsi_{,N})$ satisfies the conditions of Theorem \ref{excessdecay} with $\alpha = 1/2$, $R=r_02^N$, and $r^*_{1/2} \leq r_0$. Therefore, we may apply the mean-value property \eqref{mvp} to $\varphi_{b_i}^{N+1}$, which is $a$-harmonic on $B_{r_02^{N}}^+$ with no-flux boundary data on $\partial \hs \cap B_{r_02^N}$. Following our application of \eqref{mvp} by a use of \eqref{energyfirst} from Lemma \ref{energyvarphi} allows us to write
\begin{align}
\begin{split}
\label{intsublin6}
\left( \fint_{B_r^+}|\nabla \varphi_{b_i}^{N+1}|^2 \,dx \right)^{1/2} & \leq C_{Mean}\left( \fint_{B_{r_02^N}^+}|\nabla \varphi_{b_i}^{N+1}|^2 \,dx \right)^{1/2}\\
&\leq C_{Mean} C_1 2^d \delta_{r_02^{N+2}}^{1/3}. 
\end{split}
\end{align}
\end{proof}

\noindent \textbf{Remark}: The significance of the previous lemma is that it shows that for a specific choice of $r_0$ the relation \eqref{assumption12} holds for all $n \geq -1$ and $r \geq r_0$. In particular, this is seen through an induction on $n$: as the $\varphi_{b_i}^{n}$ for $n \in \left\{-1,0,1,2,3 \right\}$ are near-field contributions for all $r \geq r_0$, Lemma \ref{energyvarphi} provides the base case for this claim and Lemma \ref{chooser0} is the inductive step.\\

\noindent \textbf{Step 4: Passing to the limit $N \rightarrow \infty$.}\\

\noindent \textit{Proof of Theorem 1.} We split this proof into two parts: first we finish constructing the half-space-adapted corrector and then, in a second step, we obtain the half-space-adapted vector potential.\\

\noindent \textbf{Part 1:} The half-space-adapted corrector $\phi^{\mathbb{H}}_{b_i}$.\\

The remark following Lemma \ref{chooser0} implies that for any $r \geq r_0$ the inequality
\begin{align}
\begin{split}
\sum_{n=-1}^{\infty}\left( \fint_{B_r^+}|\nabla \varphi_{b_i}^n|^2 \,dx\right)^{1/2} \leq C_2 \sum_{n=n_0}^{\infty} \delta_{2^{n}}^{1/3}
\end{split}
\end{align}
holds. Furthermore, the Poincar\'{e} inequality with zero average applied in the form \eqref{averageball1poincare} with $u = \varphi_{b_i}^n$ gives that 
\begin{align}
\begin{split}
\sum_{n=-1}^{\infty}\left( \fint_{B_r^+}|\varphi_{b_i}^n |^2 \,dx\right)^{1/2} & \lesssim r^{(d+2)/2} \sum_{n=-1}^{\infty} \left(\fint_{B_r^+} |\nabla \varphi_{b_i}^n|^2 \,dx\right)^{1/2}.
\end{split}
\end{align}
Therefore, thanks to our assumption \eqref{reqwsc} on the whole-space corrector/vector potential pair, the sum $\sum_{n=-1}^{\infty} \varphi_{b_i}^n$ converges absolutely in $H^1(B_r^+)$; We call the limit $\varphi_{b_i}$.

To see that $\varphi_{b_i}$ is sublinear we again use the remark following Lemma \ref{chooser0} combined with the Poincar\'{e} inequality with zero average to obtain
\begin{align}
\label{rate1}
\begin{split}
\frac{1}{r}\left( \fint_{B_r^+}\left|\varphi_{b_i} - \fint_{B_r^+} \varphi_{b_i}\, dx\right|^2 \, dx\right)^{1/2} & \leq  C_P \left( \fint_{B_r^+} |\nabla \varphi_{b_i}|^2 \,dx \right)^{1/2}\\
& \leq C_P \sup_{N} \sum_{n=-1}^N \left( \fint_{B_r^+}|\nabla \varphi_{b_i}^n|^2 \,dx \right)^{1/2}\\
& \leq C_P C_2 \sum_{n=n_0}^{\infty}\min\left\{1, \left(\frac{2^n}{r}\right)^{d/2} \right\}\delta_{2^n}^{1/3},
\end{split}
\end{align}
which is sufficient due to \eqref{reqwsc} and the dominated convergence theorem for sums. 

To conclude, we find that the half-space-adapted corrector can be taken to be $\hp_{b_i} = \phi_{b_i} + \varphi_{b_i}$ for $i \neq d$ and $\phi_{b_d}^{\mathbb{H}} = \phi_{b_d} |_{\mathbb{H}_+^d}$. The relation \eqref{rate1} and (\ref{reqwsc}) yield the desired sublinearity property 
\begin{align}
\label{correction2}
\begin{split}
\lim_{r\rightarrow\infty} \frac{1}{r} \left( \sum_{i=1}^{d-1} \fint_{B_r^+}\left|\hp_{b_i} - \fint_{B_r^+}\hp_{b_i}\,dx\right|^2 \,dx  + \fint_{B_r}|\phi_{b_d}^{\mathbb{H}}|^2 \,dx\right)^{1/2} =0.
\end{split}
\end{align}

\noindent \textbf{Part 2:} The half-space-adapted vector potential $\hsi_{b_i}$.\\

We now pass to the limit $N \rightarrow \infty$ in the sequence $\{\psi_{b_ijk,N} \}_N$ by showing that it is a Cauchy sequence in $L^2(B_r^+)$ for all $r >r_0$.  First, we notice that $v_{b_i j, N+1}^n - v^n_{b_i j,N}$ satisfies the equation
  \begin{subequations}
    \label{equndifference}
    \begin{align}
    - \Delta (v_{b_i j, N+1}^n - v^n_{b_i j,N}) &=  \di ( \eta_n x_j a \nabla \varphi_{b_i}^{N+1}) &&\textrm{in } \hs,\\
    ~~~~ v_{b_i j, N+1}^n - v^n_{b_i j,N}  & =  0 &&\textrm{for } j \neq d \textrm{ on } \bhs,\\
    ~~~~ \partial_d (v_{b_i d, N+1}^n - v^n_{b_i d,N})  & =  0 && \textrm{on } \bhs
    \end{align}
    \end{subequations} 
in a distributional sense when $j \neq d$ and in a weak sense when tested against $H^1_{bdd}(\hs)$ when $j = d$. We consider this equation for $N \geq -2$ and adopt the notation $v^n_{b_i j, -2} = 0$ and $c_{b_i j,-2}^n=0$.

Repeating the argument from $i)$ of the proof of Lemma \ref{2.2}, gives that
\begin{align}
\begin{split}
 \label{correction1}
 & \frac{1}{r}\left(\fint_{B_r^+} |\partial_k (v^n_{b_i j, N+1} - v^n_{b_i j, N})  - e_k \cdot (c_{b_i j,N+1}^n-c_{b_ij,N}^n)|^2 \, dx\right)^{1/2} \\
   &\leq  C_3 \min \left\{1, \frac{r_02^{n+1}}{r} \right\} \left( \fint_{B_{r_02^{n+1}}^+} |\nabla (\hp_{b_i, N+1} -\hp_{b_i, N}) |^2\, dx\right)^{1/2}\\
 &  =  C_3 \min \left\{1, \frac{r_02^{n+1}}{r} \right\} \left( \fint_{B_{r_02^{n+1}}^+} |\nabla \varphi^{N+1}_{b_i} |^2\, dx\right)^{1/2}.
\end{split}
\end{align}
Summing in $n$, we find that the $v_{b_ij,N}$ from Step 3 satisfy
\begin{align}
\label{pfthm12}
\begin{split}
&\frac{1}{r} \left( \fint_{B_r^+}|\partial_k( v_{b_i j, N+1} - v_{b_i j,N})|^2 \, dx \right)^{1/2}\\
& \leq C_3 \sum_{n=-1}^{\infty}\min\left\{1, \frac{r_02^{n+1}}{r} \right\}\left(\fint_{B_{r_02^{n+1}}^+}|\nabla \varphi_{b_i}^{N+1}|^2 \,dx \right)^{1/2}.
\end{split}
\end{align}
We then sum (\ref{pfthm12}) over $N$, which as we assume that $r\geq r_0$, by \eqref{assumption12} yields 
\begin{align}
\label{pfthm13}
\begin{split}
&\frac{1}{r}\sum_{N=-2}^{\infty} \left( \fint_{B_r^+}|\partial_k( v_{b_i j, N+1} - v_{b_i j,N})|^2 \, dx \right)^{1/2}\\
& \leq C_3 C_2 \sum_{N=-2}^{\infty}\sum_{n=-1}^{\infty} \min\left\{ 1, \frac{r_02^{n+1}}{r}\right\}  \min \left\{ 1, 2^{d(N+1-n)/2}\right\}\delta_{r_02^{N+2}}^{1/3}.
\end{split}
\end{align}
To complete our argument we notice that 
\begin{align}
\label{pfthm14}
\begin{split}
& \sum_{N=-2}^{\infty}\sum_{n=-1}^{\infty} \min\{1, 2^{d(N+1-n)/2}\} \delta_{r_02^{N+2}}^{1/3}\\ 
& \leq \sum_{N=n_0}^{\infty} (N + \frac{1}{1-2^{-d/2}})\delta_{2^N}^{1/3}\\
& < \infty,
\end{split}
\end{align}
where we have used the assumption \eqref{reqwsc}. By (\ref{pfthm13}), (\ref{pfthm14}), and the definition $\psi_{b_i jk, N} = \partial_k v_{b_i j,N} - \partial_j v_{b_i k,N}$ we find that $\{\psi_{b_i jk, N}\}_N$ is a Cauchy sequence in $L^2(B_r^+)$ for all $r>0$. We may, therefore, on every half-ball $B_r^+$ pass to the limit, which we denote as $\psi_{b_i jk}$.

Also following from (\ref{pfthm13}) and (\ref{pfthm14}), this time using the dominated convergence theorem for sums, is the sublinearity property:
\begin{align}
\label{pfthm15}
\begin{split}
&\lim_{r\rightarrow \infty}\frac{1}{r}\left( \fint_{B_r^+} |\psi_{b_i jk}|^2 \, dx \right)^{1/2}\\
&\leq \lim_{r \rightarrow \infty } \frac{2}{r} \sup_{N \geq -1} \sup_{j,k \in \{1,...,d \}} \left( \fint_{B_r^+} |\partial_k v_{b_i j,N}|^2 \, dx \right)^{1/2}\\
&=0.
\end{split}
\end{align}

Recall that the proposed half-space-adapted vector potential for $ i \neq d$ is $\hsi_{b_i jk}= \sigma_{b_i jk} + \psi_{b_i jk}$. In the discussion proceeding Lemma 2.2 we concluded that for $N \geq -1$ the intermediate half-space-adapted vector potentials $\hsi_{b_i j, N}$ are distributional solutions to (\ref{partialhalfspacepotential}). Now, for any test function $\psi \in C_0^{\infty}(\hs)$  we may find $r >0$ such that $\textrm{supp}(\psi) \subset B_{r}^+$. Since the $\hsi_{b_ijk,N}$ converge to $\hsi_{b_ijk}$ in $L^2(B_{r}^+)$ and the $\hp_{b_i,N}$ converge to $\hp_{b_i}$ in $H^1(B_{r}^+)$, we find that (\ref{halfspacepotential}) is satisfied. \\

In summary, the half-space-adapted corrector/vector potential pair $(\hp, \hsi)$ that we have constructed is given by: $\hp= (\hp_{b_1},...,\hp_{b_d})$, with $\hp_{b_i}$ taken to be the limit from Part 1 when $i \in\{1,...,d-1\}$ and $\hp_{b_d}=\phi_{b_d}|_{\hs}$, and $\hsi_{jk} = (\hsi_{b_1jk}, ..., \hsi_{b_d jk})$, where $\hsi_{b_i jk}$ is the limit from Part 2 when $i \in \{1,...,d-1 \}$ and $\hsi_{b_djk} = \sigma_{b_d jk}|_{\hs}$. That the sublinearity condition (\ref{sublinear}) is satisfied can be seen from (\ref{correction2}), (\ref{pfthm15}), and the condition (\ref{reqwsc}).
\begin{flushright}
$\Box$
\end{flushright}

\section{Argument for the Excess Decay}

We first recall two basic lemmas. The first is a Caccioppoli estimate for $a$-harmonic functions with no-flux boundary data on the half-space.

\begin{customlemma}{3.1}\label{Caccioppolifornoflux} Let $a \in \Omega$, where $\Omega$ is defined by \eqref{Omega}, and $r>0$. Then for any function $u$ that is $a$-harmonic on $B_{2r}^+$ and has no-flux boundary conditions on $\fpb_{2r}$ the estimate
\begin{align}
 \label{innerregcac}
 \int_{B_r^+} |\nabla u|^2 \,dx \lesssim \frac{1}{r^2} \int_{B_{2r}^+}|u|^2 \,dx
\end{align}
\noindent holds.
\end{customlemma}

\begin{proof} Let $\eta$ denote a radial cut-off such that $\eta(x) \equiv 1$ when $|x|\leq r$, $\eta(x) \equiv 0$ when $ |x|\geq 2r$, and $|\nabla \eta (x)| \leq  \frac{2}{r}$. We test the equation 
 \begin{subequations}
 \label{aharmonic}
 \begin{align}
  ~~~~ -\di (a \nabla u) & =0 && \textrm{in } \quad B_{2r}^+,\\
  \label{nofluxbc}
  ~~~~ e_d \cdot a (\nabla u) & = 0 && \textrm{on } \quad \fpb_{2r}
 \end{align}
\end{subequations}
\noindent with $\eta^2 u$. The boundary term vanishes on $\fpb_{2r}$ due to the no-flux boundary condition (\ref{nofluxbc}) and also on $\rpb_{2r}$ due to the cut-off $\eta$. Using the uniform ellipticity and boundedness of $a$ and Young's inequality gives
\begin{align}
 \label{innerregcac1}
 \begin{split}
  \lambda \int_{B_{2r}^+} \eta^2 |\nabla u|^2 \, dx \leq \int_{B_{2r}^+} \frac{\lambda}{2}\eta^2|\nabla u|^2 + \frac{2}{\lambda} |\nabla \eta|^2 u^2 \, dx.
 \end{split}
\end{align}
\noindent To finish the argument one absorbs the first term on the right-hand side of (\ref{innerregcac1}) into the left-hand side and uses the properties of $\eta$.
\end{proof}

We will need the following facts about \textit{constant coefficient}-harmonic functions:

\begin{customlemma}{3.2}\label{innerreglemma} Let $a_{const} \in \Omega$ be constant, where $\Omega$ is defined by \eqref{Omega}, and fix $R > 0$. Let $v$ be $a_{const}$-harmonic on $B_{R}^+$ with no-flux boundary conditions on $\fpb_{R}$; i.e., $v$ solves
  \begin{subequations}
    \label{equationforv}
    \begin{align}
    ~~~~~~~~-\di (a_{const} \nabla v) &= 0  &&\textrm{in }\quad B_{R}^+,\\ 
    \label{dirichletbcv}
    ~~~~~~~~ v & = u  && \textrm{on } \quad \rpb_{R},\\
    \label{nofluxbcv}
    ~~~~~~~~ e_d \cdot a_{const} \nabla v & =  0 &&\textrm{on }\quad \fpb_{R}
    \end{align}
    \end{subequations}
\noindent for some function $u\in H^{1/2}(\partial B_r^+)$. Then for any positive $\rho \leq \frac{R}{2}$ and $r \leq \frac{R}{2}$ there exists a $\beta(d, \lambda)>0$ such that the following 
estimates hold:
\begin{subequations}
\label{innerreg}
\begin{align}
\label{innerreg1}
 \sup_{x \in B_r^+} |\nabla^n v(x)|^2 & \lesssim   \left(\frac{1}{R}\right)^{2(n-1)} \fint_{B_{R}^+}|\nabla v|^2 \, dx \quad \quad \quad \hspace{.3cm} \textrm{ for any $n \geq 1$},\\
 \label{innerreg2}
 \int_{A^{\prime}}|\nabla v|^2 \, dx& \lesssim  (R)^{1-\beta} \rho^{\beta} \int_{\rpb_{R} }|\nabla^{tan} u|^2 \, dS,\\
 \label{innerreg3}
 \textrm{ and }\sup_{x \in A^{\prime \prime}} |\nabla^n v(x)|^2 & \lesssim \left(\frac{1}{\rho}\right)^{2(n-1)} \left( \frac{R}{\rho}\right)^d \fint_{B_{R}^+}|\nabla v|^2 \, dx \quad \textrm{ for any $n \geq 1$},
\end{align}
\end{subequations}
\noindent where we have used the notation 
\begin{align}
\begin{split}
&A^{\prime} =  (B_{R}^+ \setminus B_{R - 2 \rho}^+) \cup (B_{R}^+ \cap \{\hspace{.1cm} x \hspace{.1cm} | \hspace{.1cm} x_d \leq 2 \rho \}) \textrm{ and }  \\
 &A^{\prime \prime}  = B_{R - \rho}^+ \setminus \{\hspace{.1cm} x \hspace{.1cm} | \hspace{.1cm} x_d \leq \rho  \}.
\end{split}
\end{align}

\end{customlemma}

\begin{proof}
The third estimate (\ref{innerreg3}) follows from the observation that for all $x \in A^{\prime \prime}$ we have the inner regularity estimate 
\begin{align}
 \label{innerregpf3.1}
 \sup_{y \in B_{\rho/2}(x)} |\nabla^n v(y)|^2 & \lesssim \frac{1}{\rho^{d+2(n-1)}} \int_{B_{\rho}(x)}|\nabla v|^2 \, dx.
\end{align}
\noindent This follows from an application of the Sobolev embedding and noting that all of the components of $\nabla^n v$ are $a_{const}$-harmonic in $B_{\rho}(x)$. We obtain (\ref{innerreg3}) by writing:
\begin{align}
\label{innerregpf3.2}
\begin{split}
 \sup_{x \in A^{\prime \prime}} |\nabla^n v(x)|^2 & \leq \sup_{x \in A^{\prime \prime}} \sup_{y \in B_{\rho/2}(x)}|\nabla^n v(y)|^2\\
 & \lesssim \sup_{x \in A^{\prime \prime}} \frac{1}{\rho^{d+2(n-1)}} \int_{B_{\rho}(x)}|\nabla v|^2 \, dx\\
 & \lesssim \frac{1}{\rho^{2(n-1)}} \left( \frac{R}{\rho} \right)^d \fint_{B_{R}^+}|\nabla v|^2 \, dx.
\end{split}
\end{align}

The first estimate (\ref{innerreg1}) is shown in a similar manner. In particular, we again use the Sobolev embedding and iterate the Caccioppoli inequality \eqref{innerregcac} by differentiating \eqref{equationforv}. However, this procedure only yields the inequality \eqref{innerregcac} for higher derivatives involving at most one derivative in the $e_d$ direction (as $\partial_d v$ does not satisfy \eqref{nofluxbcv}). Using a standard argument, one obtains the required estimates for higher derivatives involving multiple normal derivatives. In particular, one expresses $\partial_d^n v$ in terms of $\partial^{\beta} v$ where $|\beta| = n$ and $\beta_d = n-1$ by using the equation \eqref{equationforv} and proceeds inductively.

To show the second estimate we extend $v$ to $B_R$ through means of even reflection across $\bhs$. The extended function (which we again call $v$) is then $\tilde{a}_{const}$-harmonic on $B_{R}$ where the now \textit{not} constant coefficients are given by
\begin{align*}
(\tilde{a}_{const})_{ij}=
\begin{cases}
(a_{const})_{ij} &\text{ for }x_d \geq 0,
\\
(a_{const})_{ij}&\text{ for }x_d<0 \text{ and }i\neq d, j\neq d,
\\
-(a_{const})_{ij}&\text{ for }x_d<0 \text{ and }i=d, j\neq d,
\\
-(a_{const})_{ij}&\text{ for }x_d<0 \text{ and }i\neq d, j= d,
\\
(a_{const})_{ij}&\text{ for }x_d<0 \text{ and }i=j=d.
\end{cases}
\end{align*}
\noindent Let $\bar{v}$ be the harmonic extension of $v|_{\partial B_{R}}$ onto $B_{R}$ and notice that the estimate $\| \nabla \bar{v} \|_{L^{2/(1-\beta)}(B_{R})} \lesssim (R)^{1/2-d \beta /2}\|\nabla^{tan} v \|_{L^2(\partial B_{R})}$ holds for $\beta>0$ small enough. This can be seen by interpolating between $\| \nabla \bar{v} \|_{L^{2d/(d-1)} (B_{R})} \lesssim \|\nabla^{tan} v \|_{L^2(\partial B_{R})}$ and $\| \nabla \bar{v} \|_{L^{2}(B_{R})} \lesssim (R)^{1/2}\|\nabla^{tan} v \|_{L^2(\partial B_{R})}$. This is combined with Meyer's estimate which tells us that, also for small enough $ \beta > 0$ (where ``small enough'' depends on $d$ and $\lambda$, see \cite{Meyers}), $\| \nabla (v  - \bar{v}) \|_{L^{2/(1-\beta)}(B_{R})} \lesssim \| \nabla \bar{v} \|_{L^{2/(1-\beta)}(B_{R})}$.

The proof of (\ref{innerreg2}) then follows from an application of H\"{o}lder's inequality:
\begin{align}
\begin{split}
  \int_{A^{\prime}}|\nabla v|^2 \, dx \leq & | A^{\prime}|^{\beta} \left( \int_{B_{R}^+}|\nabla v|^{2/(1- \beta)}\, dx \right)^{1-\beta}\\
  \lesssim & (R)^{1-\beta} \rho^{\beta} \int_{\rpb_{R}}|\nabla^{tan} u|^2 \, dS.
  \end{split}
\end{align}
\end{proof}

\begin{proof}[Proof of Theorem 2]\hspace{1cm}\\

\noindent \textbf{Step 1: Main ingredient for the proof of the half-space-adapted excess-decay.}\\

First, notice that due to the linearity of the map $\xi \mapsto \hp_{\xi}$ we may express
\begin{align}
\label{excessdefn}
\exch(r) = \inf_{\xi \in \mathbb{R}^d} \displaystyle\fint_{B_r^+}\left|\nabla u - \displaystyle\sum_{i=1}^{d-1} \langle b_i , \xi \rangle (b_i + \nabla \hp_{b_i}) \right|^2 \, dx.
\end{align}
Using this form of the half-space-adapted tilt-excess, we sum up the main ingredient of the proof as:\\

\noindent \textit{Claim:} There exists $\beta>0$ and a radius $r^{\prime}>0$ such that for any radius $r$ such that $r^{\prime} \leq r \leq  R$ there exists a $\xi \in \mathbb{R}^d$ such that
\begin{align}
\label{mainexcess}
 \begin{split}
   & \fint_{B_{r}^+} \left|\nabla u - \sum_{i=1}^{d-1} \langle b_i, \xi\rangle (b_i + \nabla \hp_{b_i})\right|^2  \, dx\\
  \lesssim & \left(  \left( \frac{R}{r} \right)^{2(d+1)} \delta^{\beta/ (d+3)} + \left(\frac{r}{R}\right)^2 \right) \fint_{B_{R}^+}|\nabla u|^2 \, dx,
\end{split}
 \end{align}
\noindent where $\delta= \max(\deltah_{R},\deltah_{2r})$. \\

Our aim in Step 1 is to prove the above claim. To set-up the argument we make a couple of simplifications and introduce some definitions: 
\begin{enumerate}
 \item [i)] Notice that \eqref{mainexcess} is clear for $ r \in \left( R/4, R \right]$ with the choice $\xi = 0$. So, we may assume that $r \leq R/4$.\\
\item [ii)] We let $\Rp \in \left[ R/2, R \right)$ be a radius such that 
\begin{align}
 \label{chooseRprime}
 \int_{\rpb_{\Rp}} |\nabla^{tan} u |^2 \, dS \leq \frac{1}{R} \int_{B_{R}^+ \setminus B_{R/2}^+}|\nabla^{tan}u|^2 \, dx ,
\end{align}
\noindent which can be seen to exist by writing the second integral in polar coordinates.\\
\item[iii)] We use two smooth cut-offs: First, a one-dimensional cut-off $L(x) = L(x_d)$ that satisfies $L(x_d) \equiv 1$ if $|x_d| \leq \rho$ and $L(x_d) \equiv 0$ if $|x_d| \geq 2 \rho$. Second, a function $\eta$ that satisfies $\eta(x) \equiv 1$ if $|x| \leq \Rp - 2 \rho$ and $\eta(x)\equiv 0$ if $|x|\geq \Rp - \rho$. We assume that $0< \rho \leq r/2$ and both $|\nabla L_d|\leq 2/\rho$ and $|\nabla \eta|\leq 2 / \rho$.\\
\end{enumerate}

The core of the argument is to consider $u$ as a perturbation of $v$ satisfying (\ref{equationforv}) from Lemma \ref{innerreglemma} with the coefficients $a_{const} \in \mathbb{R}^{d \times d}$ taken to be the homogenized coefficients $a_{hom}$. In particular, $v$ is taken to satisfy
\begin{subequations}
\begin{align}
\label{equationvexcess}
-\di (a_{hom} \nabla v) &= 0  \quad \quad \quad \textrm{in }\quad B_{\Rp}^+,\\
\label{vexcessbcr}
v & = u  \quad \quad \quad \textrm{on } \quad \rpb_{\Rp},\\
\label{vexcessbcf}
e_d \cdot a_{hom} \nabla v & =  0 \quad \quad \quad \textrm{on } \quad \fpb_{\Rp}.
\end{align}
\end{subequations}
Interpreting the boundary condition \eqref{vexcessbcf} in the distributional sense, one may find a solution $v \in H^1(B_{\Rp}^+)$ to this equation using a Lax-Milgram argument.

Then notice that thanks to the remark in the proof of Lemma 2.2 we may actually interpret \eqref{vexcessbcf} in a point-wise sense. Decomposing $\nabla v = \langle b_i, \nabla v \rangle b_i $ and using that $b_i \in B$ when $i \neq d$ gives that
\begin{align}
 0 = e_d \cdot a_{hom} \nabla v(0) = e_d \cdot a_{hom} b_d \langle b_d, \nabla v(0) \rangle.
\end{align}
As $e_d \cdot a_{hom} b_d \neq 0$ this implies that $\langle b_d, \nabla v(0) \rangle =0$.

Having made this observation, we may use Young's inequality to write
\begin{align}
\label{terms}
\begin{split}
& \, \int_{B_{r}^+} \left|\nabla u - \sum_{i=1}^{d-1} \langle b_i, \nabla v(0)\rangle (b_i + \nabla \hp_{b_i})\right|^2 \, dx\\
\lesssim &   \int_{B_{r}^+} |\nabla (u - v) - (1-L) \langle b_i, \nabla v \rangle  \nabla \hp_{b_i} )|^2 \, dx\\
&  +  \int_{B_{r}^+} |\langle b_i,  \nabla v -\nabla v(0)\rangle (b_i + \nabla \hp_{b_i})|^2 \, dx\\
& + \int_{B_{r}^+} |L\langle b_i, \nabla v \rangle  \nabla \hp_{b_i}|^2 \, dx\\
 \end{split}
\end{align}
\noindent and then treat the three terms on the right-hand side separately.

We begin with the first term. Let $w= u - v -\eta (1-L) \langle b_i, \nabla v \rangle \hp_{b_i}$ denote the ansatz for the homogenization error given by two-scale expansion. Since $r \leq \Rp- 2 \rho$ we have that
\begin{align}
 \label{firstterm}
 \begin{split}
 &\int_{B_{r}^+} |\nabla(u-v) - (1-L)\langle b_i, \nabla v \rangle  \nabla \hp_{b_i})|^2 \, dx\\
 \lesssim &  \int_{B_{r}^+}|\nabla w|^2 \, dx+ \int_{B_{r}^+}|\nabla( (1-L) \langle b_i, \nabla v \rangle ) \hp_{b_i}|^2 \, dx. 
 \end{split}
\end{align}

Using that  $w \in H^1(B_{\Rp}^+)$ with $w = 0$ on $\rpb_{\Rp}$, the equations for the half-space-adapted corrector and vector potential, and  the properties of $u$ and $v$, we estimate $\|\nabla w \|_{L^2({B_{\Rp}^+})}$:
\begin{align*}
  & \int_{B_{\Rp}^+}\nabla w \cdot a \nabla w \, dx\\
  & ~~~=  \int_{B_{\Rp}^+} \nabla w \cdot a \nabla(u - v - \eta (1-L) \langle b_i, \nabla v \rangle \hp_{b_i})\, dx\\
  & ~~~=  - \int_{B_{\Rp}^+} \nabla w \cdot a \nabla (v + \eta (1-L) \langle b_i, \nabla v \rangle \hp_{b_i})\, dx\\
  & ~~~= - \int_{B_{\Rp}^+} \nabla w \cdot  ((1 - \eta(1-L)) a \nabla v + a \nabla(\eta (1-L) \langle b_i, \nabla v \rangle) \hp_{b_i})\, dx\\
  & ~~~~~~~~~~~~~~~ - \int_{B_{\Rp}} \nabla w \cdot   \eta (1-L)  \langle b_i, \nabla v \rangle  a (b_i + \nabla \hp_{b_i})\\
  & ~~~= - \int_{B_{\Rp}^+}  \nabla w \cdot   ((1 - \eta(1-L)) a \nabla v + a \nabla(\eta (1-L) \langle b_i, \nabla v \rangle) \hp_{b_i}  ) \, dx \\
  & ~~~~~~~~~~~~~~~ + \int_{B_{\Rp}^+} w \nabla (\eta (1-L) \langle b_i, \nabla v \rangle) \cdot a (b_i + \nabla \hp_{b_i}) \, dx \\
  & ~~~ = - \int_{B_{\Rp}^+}  \nabla w \cdot  (a-a_{hom}) (1 - \eta(1-L)) \nabla v + \nabla w \cdot a \nabla(\eta (1-L) \langle b_i, \nabla v \rangle) \hp_{b_i} \, dx \\
  & ~~~~~~~~~~~~~~~ + \int_{B_{\Rp}^+} w \nabla (\eta (1-L) \langle b_i, \nabla v \rangle) \cdot (a (b_i + \nabla \hp_{b_i}) - a_{hom} b_i) \, dx\\
   & ~~~= - \int_{B_{\Rp}^+} \left( \nabla w \cdot  (a-a_{hom}) (1 - \eta(1-L)) \nabla v + \partial_k w  \partial_j (\eta (1-L) \langle b_i, \nabla v \rangle) \sigma_{b_ijk} \right.\\
  & ~~~~~~~~~~~~~~~ \left. + \nabla w \cdot a \nabla(\eta (1-L) \langle b_i, \nabla v \rangle) \hp_{b_i}  \right) \vphantom{\int_{B_{\Rp}^+}} \, dx.
\end{align*}
Notice that the last step follows from the skew-symmetry of $\sigma_i$ and uses the Einstein summation convention. Also, we remark that we may test the equation \eqref{halfspacepotential} for the half-space-adapted vector potential with $w \nabla (\eta(1-L)\langle b_i, \nabla v \rangle )$ thanks to the presence of the cut-off functions; In particular, $w \nabla (\eta(1-L)\langle b_i, \nabla v \rangle ) \in H^1_0(B_{\Rp}^+)$.

After using H\"{o}lder's inequality and the uniform ellipticity and boundedness of $a$ and $a_{hom}$, this calculation (in the notation from Lemma \ref{innerreglemma}) gives that
\begin{align}
 \label{energyw1}
 \begin{split}
 \int_{B_{\Rp}^+} |\nabla w|^2 \, dx & \lesssim \int_{A^{\prime}} |\nabla v|^2 \, dx + \sup_{x \in A^{\prime \prime}}(|\nabla^2 v|^2 + \frac{1}{\rho^2}|\nabla v|^2)\\
  & ~~~~~~~~~~~~~~~~~~ \times \displaystyle\sum_{i=1}^{d-1} \int_{A^{\prime \prime}} |(\hp_{b_i} - \fint_{B_R^+} \hp_{b_i} \, dx , \hsi_{b_i})|^2 + |\hp_{b_d}, \hsi_{b_d}|^2 \, dx.  \\
\end{split}
\end{align}
\noindent Conveniently, we may also bound the second term of (\ref{firstterm}) in terms of the second term of (\ref{energyw1}). Applying (\ref{innerreg2}) and (\ref{innerreg3}) for $n=1$ and $n=2$ and using that we have chosen $\Rp$ according to (\ref{chooseRprime}), then gives that  
\begin{align}
 \label{energyw3}
 \begin{split}
  \int_{B_{r}^+} |\nabla (u-v) -  (1-L)\langle b_i, \nabla v \rangle  \nabla \hp_{b_i}|^2  \, dx \lesssim &  \left(\frac{\rho}{R}\right)^{\beta} \int_{B_{R}^+}|\nabla u|^2 \, dx\\
 & + \left(\frac{R}{\rho}\right)^{d+2} (\deltah_{R})^2 \int_{B_{\Rp}^+} |\nabla v|^2 \, dx.
\end{split}
\end{align}

We continue and bound the second term on the right-hand side of (\ref{terms}). Here, an application of (\ref{innerreg1}) for $n=2$ yields
\begin{align}
\label{termssecondterm1}
 \begin{split}
  \int_{B_{r}^+} |\langle b_i,  \nabla v -\nabla v(0)\rangle (b_i + \nabla \hp_{b_i})|^2 \, dx & \lesssim r^2 \sup_{x \in B_{r}^+} |\nabla^2 v|^2 \int_{B_{r}^+}|b_i + \nabla \hp_{b_i}|^2 \, dx\\
  & \lesssim \left( \frac{r}{R}\right)^2 \fint_{B_{\Rp}^+}|\nabla v|^2 \, dx \, \int_{B_{r}^+}|b_i + \nabla \hp_{b_i}|^2 \, dx.
 \end{split}
\end{align}
Notice that for $i = d$ the whole-space Caccioppoli estimate and for $i \neq d$ the estimate (\ref{innerregcac}) together imply that
\begin{align}
 \label{termssecondterm2}
 \begin{split}
  \fint_{B_{r}^+}|b_i + \nabla \hp_{b_i}|^2 \, dx &  \lesssim 1 + (\deltah_{2r})^2.
  \end{split}
\end{align}
The combination of \eqref{termssecondterm1} and \eqref{termssecondterm2} then gives:
\begin{align}
 \label{termssecondterm4}
 \begin{split}
  \int_{B_{r}^+} |\langle b_i,  \nabla v -\nabla v(0)\rangle (b_i + \nabla \hp_{b_i})|^2 \, dx & \lesssim  r^d (1 + (\deltah_{2r})^2) \left( \frac{r}{R}\right)^2 \fint_{B_{\Rp}^+}|\nabla v|^2 \, dx.
 \end{split}
\end{align}

Lastly, we treat the third term on the right-hand side of (\ref{terms}). An application of (\ref{innerreg1}) for $n=1$ gives

\begin{align}
\label{termsthirdterm1}
 \int_{B_r^+} |L \langle b_i, \nabla v \rangle  \nabla \hp_{b_i}|^2 \, dx \lesssim \fint_{B_{\Rp}^+} |\nabla v|^2 \, dx  \int_{{B_r^+} \cap \left\{ x_d \leq 2 \rho \right\}} |\nabla \hp_{b_i}|^2 \, dx.
\end{align}

\noindent To treat the right-hand side of (\ref{termsthirdterm1}) we modify the box-wise Caccioppoli argument used in Lemma \ref{energyvarphi}. Using the same notation (the $d$-dimensional box with center $z\in \mathbb{R}^d$ and length $l \in \mathbb{R}$ is denoted as $C_{l}(z)$), we cover $B_{r}^+ \cap \left\{ x_d \leq 2 \rho \right\}$ with boxes of width $4 \rho$ with centers taken from a set
\begin{align}
\begin{split}
 \label{centers}
 S = & \left\{ \vphantom{\sum_{z \in S}} z \in \mathbb{R}^d \hspace{.2cm} \middle| \hspace{.2cm} |B_{r}^+ \cap \left\{ x_d \leq 2 \rho \right\} \setminus \cup_{z \in S}C_{4 \rho}(z)| = 0, \quad \cup_{z \in S}C_{6 \rho}(z) \subseteq B_{2r} \right.\\
 & \left. \textrm{ and for all } x \in \mathbb{R}^d \textrm{ we have that } \sum_{z \in S} \chi_{C_{6 \rho }(z)}(x) \leq 2^d \vphantom{\mathcal{H}^d} \right\}
\end{split}
\end{align}
\noindent Then we let $\tilde{C}_{4\rho, 6 \rho, z}$ be the cut-off of $C_{4 \rho}(z)$ in the box of side length $6 \rho$ centered around it (see (\ref{defnCjcutoff}) for the definition).
 
When $i \neq d$,  for each $z \in S$, we test the half-space corrector equation (\ref{halfspacecorrector}) with $(\tilde{C}_{4 \rho, 6 \rho, z })^2 (\hp_{b_i} +b_i \cdot (x - z))$. This gives
\begin{align}
\label{cacboxsecond1}
 \begin{split}
  \int_{C_{4 \rho}(z)\cap \hs} |\nabla \hp_{b_i} +b_i |^2 \, dx & \lesssim \frac{1}{\rho^2} \int_{C_{6\rho}(z)\cap \hs} |\hp_{b_i} + b_i \cdot (x - z) - \fint_{B_{2r}^+}\hp_{b_i} \, dx|^2 \, dx,\\
 \end{split}
\end{align}
\noindent where the boundary term has vanished due to the boundary condition (\ref{boundaryconditions}). When $ i = d$ testing the whole-space corrector equation (\ref{wsceqp}) with $(\tilde{C}_{4\rho, 6 \rho, z})^2  \phi_{b_d}$ gives that
\begin{align}
\label{cacboxsecond2}
\begin{split}
 \int_{C_{4 \rho}(z)\cap \hs} |\nabla \hp_{b_d} + b_d|^2 \, dx& \lesssim \frac{1}{\rho^2} \int_{C_{6 \rho}(z)} |\phi_{b_d} + b_d \cdot (x - z)|^2 \, dx.
\end{split}
\end{align}
\noindent Summing over the $z \in S$ as in (\ref{varphienergycac3}) gives that
\begin{align}
\label{cacboxsecond3}
\begin{split}
  \int_{B_{r}^+ \cap \left\{ x_d \leq 2 \rho \right\}}|\nabla \hp_{b_i}|^2  \, dx & \lesssim r^{d-1}\rho + \left(\frac{r}{\rho}\right)^2 r^d (\deltah_{2r})^2 .
 \end{split}
\end{align}
Combining (\ref{cacboxsecond3}) with (\ref{termsthirdterm1}) then allows us to conclude:
\begin{align}
\label{termsthirdterm2}
\begin{split}
 \int_{B_{r}^+} |L\langle b_i, \nabla v \rangle  \nabla \hp_{b_i}|^2 \, dx & \lesssim  \left(r^{d-1}\rho + \left(\frac{r}{\rho}\right)^2 r^d (\deltah_{2r})^2\right) \fint_{B_{\Rp}^+} |\nabla v|^2 \, dx.
 \end{split}
\end{align}

Having treated all three terms on the right-hand side of (\ref{terms}), with the estimates (\ref{energyw3}), (\ref{termssecondterm4}), and (\ref{termsthirdterm2}), we may now write
\begin{align}
 \label{termsafter1}
 \begin{split}
   &\fint_{B_{r}^+} \left|\nabla u - \sum_{i=1}^{d-1} \langle b_i, \nabla v(0)\rangle (b_i + \nabla \hp_{b_i})\right|^2 \, dx\\
  \lesssim & \left( \left( \frac{R}{r}\right)^d   \left( \frac{R}{\rho}\right)^{d+2} (\deltah_{R})^2 + \left( \frac{r}{R} \right)^2 \left(1+ (\deltah_{2r})^2 \right) + \frac{\rho}{r} + \left( \frac{r}{\rho} \right)^2 (\deltah_{3r})^2 \right) \fint_{B_{\Rp}^+}|\nabla v|^2 \, dx\\
  & + \left(\frac{R}{r}\right)^d \left( \frac{\rho}{R} \right)^{\beta} \fint_{B_{R}^+}|\nabla u|^2 \, dx \\
  \lesssim & \left( \left( \frac{R}{r}\right)^{2(d+1)}  \left( \frac{r}{\rho}\right)^{d+2} \delta^2 + \left( \frac{r}{R} \right)^2 \left(1+ \delta^2 \right)   + \frac{\rho}{r} \right) \fint_{B_{\Rp}^+}|\nabla v|^2 \, dx\\
   & +   \left(\frac{R}{r}\right)^d  \left( \frac{\rho}{r} \right)^{\beta} \fint_{B_{R}^+}|\nabla u|^2 \, dx.
\end{split}
\end{align}
\noindent Here, we have used the notation $ \delta = \max\left\{\deltah_{R}, \deltah_{2r} \right\}$ and that $r \leq R/4$ and $\rho \leq r/2$. To post-process this estimate we do two things: derive an apriori estimate for $\|\nabla v\|_{L^2(B_{\Rp}^+)}$ and choose a specific width $\rho$ for the boundary layer introduced by the cut-offs $\eta$ and $L$.

The apriori estimate for $\nabla v$ follows from the equation satisfied by the difference $v-u$: 
  \begin{subequations}
    \label{equationfordifference}
    \begin{align}
    ~~~~~~~~-\di (a_{hom} \nabla (v-u)) &= \di a_{hom} \nabla u  &&\textrm{in }\quad B_{\Rp}^+,
    \\
    ~~~~~~~~ v-u & = 0  && \textrm{on } \quad \rpb_{\Rp},
    \\
    ~~~~~~~~ e_d \cdot a_{hom} \nabla (v-u) & =  -e_d \cdot a_{hom} \nabla u  &&\textrm{on }\quad \fpb_{\Rp}.
    \end{align}
    \end{subequations}
Testing (\ref{equationfordifference}) with $v-u$ and using H\"{o}lder's inequality then yields that
\begin{align}
 \int_{B_{\Rp}^+} |\nabla (v-u)|^2 \, dx \lesssim \int_{B_{\Rp}^+}|\nabla u|^2,
\end{align}
which by Young's inequality gives
\begin{align}
\label{aprioriestimate}
 \int_{B_{\Rp}^+} |\nabla v|^2 \, dx \lesssim \int_{B_{\Rp}^+}|\nabla u|^2.
\end{align}

We then turn to choosing the width $\rho$. Recall that the only assumption on $\rho$ was that $ \rho \in (0, r/2]$. By varying $\rho$ within this interval we may obtain $\rho/r =s$ for any $s \in (0, 1/4]$. We set $\rho$ to satisfy $\rho/r = \min \left\{ 1/4, \delta^{2/(d+3)} \right\}$.

These observations allow us to, for sufficiently large $r$ and $R$, re-write (\ref{termsafter1}) as
\begin{align}
 \begin{split}
   & \fint_{B_{r}^+} \left|\nabla u - \sum_{i=1}^{d-1} \langle b_i, \nabla v(0)\rangle (b_i + \nabla \hp_{b_i})\right|^2  \, dx\\
  \lesssim & \left(  \left( \frac{R}{r} \right)^{2(d+1)} \delta^{2 / (d+3)} + \left(\frac{r}{R}\right)^2 + \left( \frac{R}{r}\right)^d\delta^{2 \beta/(d+3)} \right) \fint_{B_{R}^+}|\nabla u|^2 \, dx.
\end{split}
 \end{align}
Notice that here ``sufficiently large $r$ and $R$'' means $R \geq r \geq r^{\prime}$ for the minimal radius $r^{\prime}>0$ guaranteeing that $\delta \leq 1$. Using that $0< \beta <1$ and $R/r \geq 1$ then yields \eqref{mainexcess}.\\

\noindent \textbf{Step 2: Proof of the half-space-adapted excess-decay.}\\

We may then apply the claim from the first step to any two radii $\tilde{r}$ and $\tilde{R}$ such that $r^{\prime} \leq \tilde{r} \leq \tilde{R}\leq R$.  Notice that thanks to (\ref{halfspacecorrector}) the function $\tilde{u}_{c} = u - \sum_{i=1}^{d-1} \langle b_i, c \rangle (b_i + \hp_{b_i})$ is $a$-harmonic with no-flux boundary conditions on $\fpb_{R}$ for any $c \in \mathbb{R}^d$. Applying \eqref{mainexcess} to these functions and taking the infimum over $c \in \mathbb{R}^d$ allows us to rephrase \eqref{mainexcess} in terms of the half-space-adapted tilt-excess:
\begin{align}
 \label{excess1}
 \exch(\tilde{r}) \lesssim \left(\theta^{-2(d+1)} \delta^{\beta/(d+3)} + \theta^{2} \right) \exch(\tilde{R}),
\end{align}
where we have used the notation $\theta = \tilde{r}/\tilde{R}$.

Thanks to condition (\ref{smallness}) and $\alpha <1$ we may choose $\theta$ and $C_{\alpha}(d, \lambda)$ such that 
\begin{align}
 \label{iterate1}
 \theta^{-2(d+1)}\delta^{\beta/(d+3)}+ \theta^{2} \leq \theta^{2 \alpha}.
\end{align}
is satisfied above some minimal radius $r^{\prime \prime} \geq r^{\prime} >0$. Making these choices, we obtain that
\begin{align}
 \exch(\theta \tilde{R} ) \lesssim \theta^{2 \alpha} \exch(\tilde{R}),
\end{align}
whenever $r^{\prime \prime} \leq \theta \tilde{R}$.

Then, for $r^{\prime \prime} \leq r \leq R$, letting $n= \lfloor \log_{\theta}(r/R)\rfloor$ we find that
\begin{align}
 \label{iterate2}
\exch(r) \leq \theta^{-d} \exch(\theta^n R) \lesssim \theta^{2 n \alpha - d} \exch(R) \lesssim \theta^{-(d+2\alpha)}\left( \frac{r}{R} \right)^{2\alpha} \exch(R).
\end{align}
This finishes the argument for the excess-decay.\\

\noindent \textbf{Step 3: Proof of the coercivity of the tilt-excess functional.}\\

We must show that 
\begin{align}
\label{Coercivity}
\fint_{B_r^+} |\nabla u -( \tilde{b} + \nabla \phi_{\tilde{b}}) |^2 \,dx \rightarrow \infty \quad \textrm{ as } \quad |\tilde{b}| \rightarrow \infty.
\end{align}
By the triangle inequality in $L^2(B_r^+)$ it suffices to prove that
\begin{align}
\label{coercive1}
\fint_{B_r^+} | \tilde{b}+ \nabla \hp_{\tilde{b}}|^2 \, dx \geq \left( \frac{1}{16}\right)^{d+1}|\tilde{b}|^2.
\end{align}

To show this we insert a smooth cut-off function $\eta$ into the left-hand side of (\ref{coercive1}), where $\eta = 1$ on $B_{r/2}^+ \cap \left\{ x_d \geq r/4 \right\}$, $\eta = 0$ outside of $B_r^+$, $0 \leq  \eta \leq 1$, and $|\nabla \eta| \leq  12/r$. It is clear that
\begin{align}
\label{coercive2}
\fint_{B_r^+} |\tilde{b}+ \nabla \hp_{\tilde{b}}|^2 \, dx \geq  \fint_{B_r^+} \eta^2 |\tilde{b}+ \nabla \hp_{\tilde{b}}|^2 \, dx.
\end{align}
\noindent Jensen's inequality and an integration by parts, in which the boundary term cancels due to the cut-off $\eta$, then yield
\begin{align}
\label{coercive3}
\begin{split}
\fint_{B_r^+} \eta^2| \tilde{b}+ \nabla \hp_{\tilde{b}}|^2 \, dx \geq & \left( \fint_{B_r^+} \eta  \,dx \right)^2 \left| \fint_{B_r^+} \frac{\eta}{\fint_{B_r^+}\eta \,dx}(\tilde{b}+ \nabla \hp_{\tilde{b}}) \, dx \right|^2\\
\geq & \left( \fint_{B_r^+} \eta  \,dx \right)^2 \left| \tilde{b} +\frac{1}{\fint_{B_r^+}\eta \,dx} \fint_{B_r^+}\eta \nabla \hp_{\tilde{b}}  \, dx\right|^2\\
= & \left( \fint_{B_r^+}  \eta \,dx \right)^2 \left|  \tilde{b} -\frac{1}{\fint_{B_r^+}\eta \,dx} \fint_{B_r^+}\nabla \eta \left( \hp_{\tilde{b}} - \fint_{B_r^+} \hp_{\tilde{b}}\, dx \right)\, dx \right|^2.
\end{split}
\end{align}
\noindent Notice that
\begin{align}
\label{coercive4}
 \left(\frac{1}{4}\right)^d \leq \fint_{B_r^+}\eta \,dx ,
\end{align}
\noindent which, along with an application of H\"{o}lder's inequality, implies
\begin{align}
\label{coercive5}
\frac{1}{\fint_{B_r^+}\eta \,dx} \left| \fint_{B_r^+}\nabla \eta  \left( \hp_{\tilde{b}}- \fint_{B_r^+} \hp_{\tilde{b}}\, dx \right) \, dx \right| \leq 4^{d+2} |\tilde{b}| \deltah_r.
\end{align}
\noindent By (\ref{smallness}) we can choose $C_{\alpha}(d, \lambda)$ large enough such that $ 4^{d+1} |\tilde{b}| \deltah_r \leq |\tilde{b}|/2 $ for all $r \geq r^{\prime \prime \prime}$. Combining this with (\ref{coercive5}) and (\ref{coercive3}), we conclude (\ref{coercive1}).\\

\noindent \textbf{Remark}: The minimal radius $r^*_{\alpha}>0$ from the statement of the theorem is then chosen to be $r^*_{\alpha} = \max(r^{\prime}, r^{\prime \prime}, r^{\prime \prime \prime})$.\\

\noindent \textbf{Step 4: Proof of the mean-value property.}\\

For any radius $r \in [ r^*_{1/2},R]$ we let $\tilde{b}_{r}\in B$ satisfy
\begin{align}
\label{infimumdefn}
\exch(r) = \fint_{B_{r}^+} |\nabla u - (\tilde{b}_{r} + \nabla \hp_{\tilde{b}_{r}})|^2 \, dx.
\end{align} 
It then holds that
\begin{align}
\label{mvp1}
\begin{split}
\fint_{B_r^+}|\nabla u|^2 \,dx \lesssim & \exch(r) + |\tilde{b}_r|^2 \\
\lesssim & \exch(R) +|\tilde{b}_r|^2\\
\lesssim & \fint_{B_R^+}|\nabla u |^2 \,dx +|\tilde{b}_R|^2+|\tilde{b}_r - \tilde{b}_R|^2.
\end{split}
\end{align}
\noindent Here, the first inequality follows from Young's inequality and \eqref{termssecondterm2}, the second uses the excess-decay from Step 2, and the third is obtained like the first. 

We must bound $|\tilde{b}_R|^2$ and $|\tilde{b}_r - \tilde{b}_R|^2$. The first bound is a simple consequence of \eqref{coercive1}, the definition of the half-space-adapted tilt-excess, and Young's inequality:
\begin{align}
\label{mvp2}
|\tilde{b}_R|^2 \lesssim \fint_{B_R^+}|\tilde{b}_R+ \nabla \hp_{\tilde{b}_R}|^2 \, dx \lesssim \exch(R) + \fint_{B_R^+}|\nabla u|^2 \, dx \lesssim \fint_{B_R^+}|\nabla u|^2 \, dx.
\end{align}
\noindent To obtain an estimate for the difference $|\tilde{b}_r-\tilde{b}_R|^2$ we first notice if $R-r \leq R/2$ then the coercivity property (\ref{coercive1}), the excess-decay,  and Young's inequality give
\begin{align}
\label{mvp3}
\begin{split}
|\tilde{b}_{r}- \tilde{b}_R|^2  \lesssim & \fint_{B_{r}^+}|\tilde{b}_{r} - \tilde{b}_{R} + (\nabla \hp_{\tilde{b}_{r}} - \nabla \hp_{\tilde{b}_R})|^2 \, dx\\
\lesssim & \exch(r) + \exch(R)\\
\lesssim & \exch(R).
\end{split}
\end{align}
Notice that the condition that $r \in [R/2, R]$ is used for the second inequality.

To finish, we iterate \eqref{mvp3}. Dropping the assumption that $r \in [R/2, R]$, let $n= \lfloor \log_{1/2} (r/R)\rfloor$. The excess-decay for $\alpha =1/2$ then gives that
\begin{align}
\label{mvp4}
\begin{split}
|\tilde{b}_r-\tilde{b}_R|^2 \leq  & \left( |\tilde{b}_r - \tilde{b}_{R2^{-n}}|+\sum_{m=1}^n|\tilde{b}_{R2^{-m}} - \tilde{b}_{R2^{-(m-1)}}|\right)^2\\
\lesssim & \left( \sum_{m=0}^n(\exch(R2^{-m}))^{1/2}\right)^2\\
\lesssim & \left( \sum_{m=0}^n2^{-m/2}(\exch(R))^{1/2}\right)^2\\
\lesssim & \hspace{.25cm} \exch(R) \hspace{.25cm} \lesssim \fint_{B_R^+} |\nabla u |^2 \,dx.
\end{split}
\end{align} 
\noindent The mean-value property then follows from (\ref{mvp1}), (\ref{mvp2}), and (\ref{mvp4}).

\end{proof}

\begin{proof}[Proof of Corollary \ref{Liouville}]  With Lemma \ref{Caccioppolifornoflux} the assumption (\ref{subquad}) of subquadratic growth can be processed to yield
\begin{align}
\label{Liouville1}
\lim_{r \rightarrow \infty} \frac{1}{r^{2\alpha} }\fint_{B_r^+}|\nabla u |^2 \, dx = 0.
\end{align}
\noindent By the definition of the half-space-adapted excess this implies that 
\begin{align}
\label{Liouville2}
\lim_{r \rightarrow \infty} \frac{1}{r^{2  \alpha}}\exch(r) =0.
\end{align}
\noindent Our condition on the whole-space corrector/ vector potential pair guarantees that the excess-decay (\ref{ExcessDecay}) holds above some minimal radius $r^*_{\alpha}>0$. So, for all $\tilde{r}>r^*_{\alpha}>0$ we have that 
\begin{align}
\label{Liouville3}
\exch(\tilde{r}) \leq \left( \frac{\tilde{r}}{r}\right)^{2\alpha} \exch(r)
\end{align}
\noindent for any $r >\tilde{r}$. Due to (\ref{Liouville2}) this implies that $\exch(\tilde{r})=0$ for all $\tilde{r} \geq r^*_{\alpha}$. Since the infimum in the definition of the half-space-adapted tilt-excess is attained,
this implies that
\begin{align}
\label{Liouville4}
 u = \tilde{b}_{\tilde{r}} \cdot x + \hp_{\tilde{b}_{\tilde{r}}} + c_{\tilde{r}} \textrm{ on } B_{\tilde{r}}^+ 
\end{align}
 \noindent for some constants $\tilde{b}_{\tilde{r}}\in \mathbb{R}^d$ and $ c_{\tilde{r}}\in\mathbb{R}$. Observe that by \eqref{mvp3} the $\tilde{b}_{\tilde{r}}$ do not depend on $\tilde{r}$.
\end{proof}

\bibliographystyle{plain}
\bibliography{stochastic_homogenization}

\begin{thebibliography}{10}

\bibitem{AleksanyanShahgholianSjolin}
H.~Aleksanyan, H.~Shahgholian, and P.~Sj{\"o}lin.
\newblock Applications of {F}ourier analysis in homogenization of the
  {D}irichlet problem: {$L^p$} estimates.
\newblock {\em Arch. Ration. Mech. Anal.}, 215(1):65--87, 2015.

\bibitem{AllaireAmar}
G.~Allaire and M.~Amar.
\newblock Boundary layer tails in periodic homogenization.
\newblock {\em ESAIM Control Optim. Calc. Var.}, 4:209--243, 1999.

\bibitem{ArmstrongKuusiMourratComplete}
S.~Armstrong, T.~Kuusi, and J.-C. Mourrat.
\newblock The additive structure of elliptic homogenization.
\newblock {\em arXiv Preprint}, 2016.
\newblock arXiv:1602.00512.

\bibitem{ArmstrongKuusiMourratPrange}
S.~Armstrong, T.~Kuusi, J.-C. Mourrat, and C.~Prange.
\newblock Quantitative analysis of boundary layers in periodic homogenization.
\newblock {\em arXiv Preprint}, 2016.
\newblock arXiv:1607.06716.

\bibitem{ArmstrongMourrat}
S.~Armstrong and J.-C. Mourrat.
\newblock Lipschitz regularity for elliptic equations with random coefficients.
\newblock {\em Arch. Ration. Mech. Anal.}, 219(1):255--348, 2016.

\bibitem{ArmstrongShen}
S.~Armstrong and Z.~Shen.
\newblock Lipschitz estimates in almost-periodic homogenization.
\newblock {\em Comm. Pure Appl. Math.}, 69(10):1882--1923, 2016.

\bibitem{ArmstrongSmart}
S.~Armstrong and C.~Smart.
\newblock Quantitative stochastic homogenization of convex integral
  functionals.
\newblock {\em Ann. Sci. \'Ec. Norm. Sup\'er}, 48:423--481, 2016.

\bibitem{AvellanedaLinCPAM}
M.~Avellaneda and F.-H. Lin.
\newblock Compactness methods in the theory of homogenization.
\newblock {\em Comm. Pure Appl. Math.}, 40(6):803--847, 1987.

\bibitem{AvellanedaLinJMPA}
M.~Avellaneda and F.-H. Lin.
\newblock Homogenization of poisson's kernel and applications to boundary
  control.
\newblock {\em J. Math. Pure Appl.}, 68:1--29, 1989.

\bibitem{BenjaminiCopinKozmaYadin}
I.~Benjamini, H.~Duminil-Copin, G.~Kozma, and A.~Yadin.
\newblock Disorder, entropy and harmonic functions.
\newblock {\em Ann. Probab.}, 43(5):2332--2373, 2015.

\bibitem{BensoussanLionsPapanicolaou}
A.~Bensoussan, J.-L. Lions, and G.~Papanicolaou.
\newblock {\em Asymptotic analysis for periodic structures}.
\newblock AMS Chelsea Publishing, 1978.

\bibitem{CorrectorEstimatesSlowDecorrelation}
J.~Fischer and F.~Otto.
\newblock Sublinear growth of the corrector in stochastic homogenization:
  Optimal stochastic estimates for slowly decaying correlations.
\newblock {\em arXiv Preprint}, 2015.
\newblock arXiv:1508.00025.

\bibitem{FischerOtto}
J.~Fischer and F.~Otto.
\newblock A higher-order large-scale regularity theory for random elliptic
  operators.
\newblock {\em Comm. Partial Differential Equations}, 41(7):1108--1148, 2016.

\bibitem{FischerRaithel}
J.~Fischer and C.~Raithel.
\newblock Liouville principles and a large-scale regularity theory for random
  elliptic opertaors on the half-space.
\newblock {\em Preprint}, 2016.
\newblock arXiv:1604:02717v2.

\bibitem{GerardVaretMasmoudi2}
D.~G\'{e}rard-Varet and N.~Masmoudi.
\newblock Homogenization in polygonal domains.
\newblock {\em J. Eur. Math. Soc. (JEMS)}, 13(5):1477--1503, 2011.

\bibitem{GerardVaretMasmoudi}
D.~G{\'e}rard-Varet and N.~Masmoudi.
\newblock Homogenization and boundary layers.
\newblock {\em Acta Math.}, 209(1):133--178, 2012.

\bibitem{GiaquintaMartinazzi}
M.~Giaquinta and L.~Martinazzi.
\newblock {\em An introduction to the regularity theory for elliptic systems,
  harmonic maps and minimal graphs}, volume~11 of {\em Appunti. Scuola Normale
  Superiore di Pisa (Nuova Serie) [Lecture Notes. Scuola Normale Superiore di
  Pisa (New Series)]}.
\newblock Edizioni della Normale, Pisa, second edition, 2012.

\bibitem{GloriaNeukammOtto}
A.~Gloria, S.~Neukamm, and F.~Otto.
\newblock A regularity theory for random elliptic operators.
\newblock {\em arXiv Preprint}, 2014.
\newblock arXiv:1409.2678.

\bibitem{GloriaOttoNew}
A.~Gloria and F.~Otto.
\newblock The corrector in stochastic homogenization: optimal rates, stochastic
  integrability, and fluctuations.
\newblock {\em arXiv Preprint}, 2015.
\newblock arXiv:1510.08290.

\bibitem{Homogbook}
V.~Jikov, S.~Kozlov, and O.~Oleinik.
\newblock {\em Homogenization of Differential Operators and Integral
  Functionals}.
\newblock Springer-Verlag, Berlin, 1994.

\bibitem{KenigLinShenNeumann}
C.~Kenig, F.-H. Lin, and Z.~Shen.
\newblock Homogenization of elliptic systems with {N}eumann boundary
  conditions.
\newblock {\em J. Amer. Math. Soc.}, 26(4):901--937, 2013.

\bibitem{KenigLinShenGreen}
C.~Kenig, F.-H. Lin, and Z.~Shen.
\newblock Periodic homogenization of {G}reen and {N}eumann functions.
\newblock {\em Comm. Pure Appl. Math.}, 67(8):1219--1262, 2014.

\bibitem{MarahrensOtto}
D.~Marahrens and F.~Otto.
\newblock Annealed estimates on the {G}reen function.
\newblock {\em Probab. Theory Related Fields}, 163(3):527--573, 2014.

\bibitem{Meyers}
N.~Meyers.
\newblock An {$L^p$}-estimate for the gradient of solutions of second order
  elliptic divergence equations.
\newblock {\em Ann. Scuola Norm. Sup. Pisa}, 17:189--206, 1963.

\bibitem{PiccininiSpagnolo}
L.~Piccinini and S.~Spagnolo.
\newblock On the {H}\"older continuity of solutions of second order elliptic
  equations in two variables.
\newblock {\em Ann. Scuola Norm. Sup. Pisa (3)}, 26:391--402, 1972.

\bibitem{ShenZhuge}
Z.~Shen and J.~Zhuge.
\newblock Boundary layers in periodic homogenization of {N}eumann problems.
\newblock {\em arXiv Preprint}, 2016.
\newblock arXiv:1610.05273.

\end{thebibliography}

\end{document}